\documentclass[a4paper,12pt]{amsart}

\usepackage[top=3cm,left=2.5cm,right=2.5cm,bottom=3cm]{geometry}
\usepackage{relsize}
\usepackage{lineno}
\usepackage{hyperref}
\usepackage{amsmath,amsthm,pb-diagram,amssymb,comment}
\usepackage{amsfonts,graphicx,color}
\usepackage{enumitem, fancyhdr, dsfont}
\usepackage{thmtools,cleveref}
\usepackage[normalem]{ulem}
\usepackage{stmaryrd}
\usepackage{textcomp}
\usepackage{contour}
\usepackage{mathbbol}
\usepackage{fontawesome5}
\usepackage{tikz,tikz-cd,float}
\usepackage{chemformula}
\usepackage{amsthm,thmtools,xcolor} 
\newlist{steps}{enumerate}{1}
\setlist[steps, 1]{label = Step \arabic*:}
\usepackage{multicol}
\hypersetup{
    colorlinks=true,
    linkcolor=dodger,
    filecolor=dodger,
    urlcolor=dodger,
    citecolor=dodger,
}

\setlist{topsep=0ex,itemsep=1ex}

\newcommand{\dom}{\mathrm{dom}}

\newcommand{\thzfc}{\mathsf{ZFC}}
\newcommand{\thgch}{\mathsf{GCH}}

 \newcommand{\Ed}{\mathsf{Ed}}
\newcommand{\Mbf}{\mathsf{M}}

\newcommand{\Mg}{\mathrm{Mg}}
\newcommand{\Cn}{\mathrm{Cn}}

\newcommand{\Bwf}{\mathcal{B}}

\newcommand{\Ewf}{\mathcal{E}}

\newcommand{\Hwf}{\mathcal{H}}
\newcommand{\Iwf}{\mathcal{I}}
\newcommand{\Jwf}{\mathcal{J}}

\newcommand{\Mwf}{\mathcal{M}}
\newcommand{\Nwf}{\mathcal{N}}

\newcommand{\Swf}{\mathcal{S}}

\newcommand{\bfrak}{\mathfrak{b}}
\newcommand{\cfrak}{\mathfrak{c}}
\newcommand{\dfrak}{\mathfrak{d}}
\newcommand{\efrak}{\mathfrak{e}}

\newcommand{\menos}{\smallsetminus}

\newcommand{\pts}{\mathcal{P}}

\newcommand{\frestr}{{\upharpoonright}}

\newcommand{\add}{\mathrm{add}}
\newcommand{\cov}{\mathrm{cov}}

\newcommand{\non}{\mathrm{non}}
\newcommand{\cof}{\mathrm{cof}}
\newcommand{\limdir}{\operatorname{lim\,dir}}

\newcommand{\Bor}{\mathds{B}}
\newcommand{\Cor}{\mathds{C}}
\newcommand{\Dor}{\mathds{D}}
\newcommand{\Eor}{\mathds{E}}

\newcommand{\Loc}{\mathds{LOC}}

\newcommand{\Por}{\mathds{P}}
\newcommand{\Pbb}{\mathds{P}}
\newcommand{\Qor}{\mathds{Q}}
\newcommand{\Ior}{\mathds{I}}

\newcommand{\Pre}{\Por\mathbb{r}}

\newcommand{\Qnm}{\dot\Qor}

\newcommand{\cf}{\mathrm{cf}}

\newcommand{\lh}{\ell g}

\newcommand{\la}{\langle}
\newcommand{\ra}{\rangle}

\newcommand{\Esf}{\mathsf{E}}

\newcommand{\Fr}{\mathsf{Fr}}

\newcommand{\D}{\mathsf{D}}
\newcommand{\Cbf}{\mathsf{C}}
\newcommand{\Rbf}{\mathsf{R}}

\newcommand{\Lc}{\mathsf{Lc}}
\newcommand{\Lb}{\mathbf{Lb}}

\newcommand{\vfa}{\mathfrak{v}}

\newcommand{\Predictors}{\Sigma}
\newcommand{\incrbaire}{{}^{\uparrow\omega}\omega}

\newcommand{\leqT}{\preceq_{\mathrm{T}}}
\newcommand{\eqT}{\cong_{\mathrm{T}}}

\newcommand{\gen}{\mathrm{gen}}

\newcommand{\const}{\mathrm{const}}
\newcommand{\pr}{\mathbf{pc}}

\renewcommand{\colon}{\nobreak\mskip2mu\mathpunct{}\nonscript
  \mkern-\thinmuskip{:}\allowbreak\mskip6muplus1mu\relax}
\newcommand{\set}[2]{\{#1\colon\,#2\}}
\newcommand{\largeset}[2]{\left\{#1\colon\,#2\right\}}
\newcommand{\seq}[2]{\la #1\colon\,#2\ra}

\newcommand{\baire}{{}^{\omega}\omega}
\newcommand{\baireinc}{{}^{\uparrow\omega}\omega}
\newcommand{\cantor}{{}^{\omega}2}

\newcommand\subsetdot{\mathrel{\ooalign{$\subset$\cr
  \hidewidth\hbox{$\cdot\mkern3mu$}\cr}}}

\newcommand{\vfrak}{\mathfrak{v}}
\newcommand{\sqsubm}{\sqsubset^{\rm m}}

\contourlength{0.8pt}
\newcommand{\bsp}{\allowbreak\ }
\newcommand{\comma}{\mathord,\bsp}
\newcommand{\setand}{\bsp\text{and}\bsp}
\newcommand{\setsth}[1]{\bsp\text{#1}\bsp}
\newcommand{\leqtr}{\triangleleft}
\newcommand{\nleqtr}{\ntriangleleft}
\newcommand{\nsqsubset}{\not\sqsubset}

\newcommand{\Cv}{\mathsf{Cv}}

\newcommand{\imp}{\mathrel{\mbox{$\Rightarrow$}}}

\newcommand{\lqq}{\textup{``}}
\newcommand{\rqq}{\textup{''}}

\contourlength{0.8pt}

\newcount\skewfactor
\def\mathunderaccent#1#2 {\let\theaccent#1\skewfactor#2
\mathpalette\putaccentunder}
\def\putaccentunder#1#2{\oalign{$#1#2$\crcr\hidewidth
\vbox to.2ex{\hbox{$#1\skew\skewfactor\theaccent{}$}\vss}\hidewidth}}

\newcommand{\checkedbox}{\makebox[0pt][l]{$\square$}\raisebox{.15ex}{\hspace{0.1em}$\checkmark$}}

\newenvironment{varproof}
 {\renewcommand{\qedsymbol}{\ensuremath{\dodger{\checkedbox}}}\proof}
 {\endproof}

\definecolor{ao(english)}{rgb}{0.0, 0.5, 0.0}
	\definecolor{ultramarineblue}{rgb}{0.25, 0.4, 0.96}
\definecolor{cornellred}{rgb}{0.7, 0.11, 0.11}
\definecolor{cobalt}{rgb}{0.0, 0.28, 0.67}
\definecolor{bleudefrance}{rgb}{0.19, 0.55, 0.91}
\definecolor{darkblue}{rgb}{0.0, 0.0, 0.55}
\definecolor{ferrarired}{rgb}{1.0, 0.11, 0.0}
\definecolor{brandeisblue}{rgb}{0.0, 0.44, 1.0}
\definecolor{azure(colorwheel)}{rgb}{0.0, 0.5, 1.0}
\definecolor{aqua}{rgb}{0.0, 1.0, 1.0}
\definecolor{aguamarina}{cmyk}{0.85,0,0.33,0}
\definecolor{cafe}{cmyk}{0,0.81,1,0.60}
\definecolor{canela}{cmyk}{0.14,0.42,0.56,0}
\definecolor{darkgray}{cmyk}{0,0,0,0.50}
\definecolor{emerald}{cmyk}{0.91,0,0.88,0.12}
\definecolor{fresa}{cmyk}{0,1,0.50,0}
\definecolor{gold}{cmyk}{0,0.10,0.84,0}
\definecolor{lightgray}{cmyk}{0,0,0,0.30}
\definecolor{marron}{cmyk}{0,0.72,1,0.45}
\definecolor{melon}{cmyk}{0,0.29,0.84,0}
\definecolor{ladri}{cmyk}{0,0.77,0.87,0}
\definecolor{olive}{cmyk}{0.64,0,0.95,0.40}
\definecolor{orange}{cmyk}{0,0.42,1,0}
\definecolor{peach}{cmyk}{0,0.46,0.50,0}
\definecolor{pink}{cmyk}{0,0.10,0.10,0}
\definecolor{orange}{cmyk}{0,0.42,1,0}
\definecolor{pine}{cmyk}{0.92,0,0.59,0.25}
\definecolor{purple}{cmyk}{0.45,0.86,0,0}
\definecolor{violet}{cmyk}{0.07,0.90,0,0.34}
\definecolor{craneorange}{RGB}{252,187,6}
\definecolor{red(ncs)}{rgb}{0.77, 0.01, 0.2}

\newcommand{\dodger}[1]{{\color{dodger}#1}}

\newcommand{\cyan}[1]{{\color{cyan}#1}}

\definecolor{aguamarina}{cmyk}{0.85,0,0.33,0}
\definecolor{cafe}{cmyk}{0,0.81,1,0.60}
\definecolor{canela}{cmyk}{0.14,0.42,0.56,0}
\definecolor{darkgray}{cmyk}{0,0,0,0.50}
\definecolor{emerald}{cmyk}{0.91,0,0.88,0.12}
\definecolor{fresa}{cmyk}{0,1,0.50,0}
\definecolor{gold}{cmyk}{0,0.10,0.84,0}
\definecolor{lightgray}{cmyk}{0,0,0,0.30}
\definecolor{marron}{cmyk}{0,0.72,1,0.45}
\definecolor{melon}{cmyk}{0,0.29,0.84,0}
\definecolor{ladri}{cmyk}{0,0.77,0.87,0}
\definecolor{olive}{cmyk}{0.64,0,0.95,0.40}
\definecolor{orange}{cmyk}{0,0.42,1,0}
\definecolor{peach}{cmyk}{0,0.46,0.50,0}
\definecolor{pink}{cmyk}{0,0.10,0.10,0}
\definecolor{orange}{cmyk}{0,0.42,1,0}
\definecolor{pine}{cmyk}{0.92,0,0.59,0.25}
\definecolor{purple}{cmyk}{0.45,0.86,0,0}
\definecolor{violet}{cmyk}{0.07,0.90,0,0.34}

\DeclareSymbolFont{extraup}{U}{zavm}{m}{n}
\DeclareMathSymbol{\varheart}{\mathalpha}{extraup}{86}
\DeclareMathSymbol{\vardiamond}{\mathalpha}{extraup}{87}

\definecolor{dodger}{rgb}{0.0,0.5,1.0}

\definecolor{amber}{rgb}{1.0,0.49,0.0}

\definecolor{ogreen}{RGB}{107,142,35}

\title[Constant evasion, constant prediction, and Cicho\'n's maximum]
{Adding the constant evasion and constant prediction numbers to Cicho\'n's maximum}

\author{Miguel A. Cardona}
\address[Miguel A. Cardona]{Einstein Institute of Mathematics\\
Edmond J. Safra Campus, Givat Ram\\
The Hebrew University of Jerusalem\\
Jerusalem, 91904, Israel}
\email{\href{mailto:miguel.cardona@mail.huji.ac.il}{miguel.cardona@mail.huji.ac.il}}
\urladdr{\url{https://sites.google.com/view/miacardonamo/home-page}}

\author{Miroslav Repick\'y}
\address[Miroslav Repick\'y]{Mathematical Institute\\
Slovak Academy of Sciences\\
Gre\v{s}\'akova~6\\
040\,01 Ko\v{s}ice\\
Slovak Republic}
\email{\href{mailto:repicky@saske.sk}{repicky@saske.sk}}

\author[S. Shelah]{Saharon Shelah}
\address[Saharon Shelah]{Einstein Institute of Mathematics,
Edmond J. Safra Campus, Givat Ram\\
The Hebrew University of Jerusalem, \\
Jerusalem, 91904, Israel; and
Department of Mathematics\\
Rutgers University\\
Piscataway, NJ 08854-8019, USA}

\email{\href{mailto:shelah@math.huji.ac.il}{shelah@math.huji.ac.il}}
\urladdr{\url{https://shelah.logic.at}}

\thanks{The first and third author would like to thank the Israel Science Foundation for partially supporting this research by grant 2320/23 (2023-2027); and the second author was supported by the grant VEGA 2/0104/24 of the Slovak Grant Agency VEGA and by the Slovak Research and Development Agency under the Contract no.\ APVV-20-0045}



\definecolor{sub0}{RGB}{29,32,137}
\definecolor{sub1}{RGB}{1,71,157}
\definecolor{sub2}{RGB}{1,104,183}
\definecolor{sub3}{RGB}{0,160,234}
\definecolor{sug}{RGB}{0,154,68}
\definecolor{suy}{RGB}{208,219,1}

\newcommand{\subiii}[1]{{\color{sub3}#1}}

\DeclareUnicodeCharacter{0301}{\'{e}}

\begin{document}

\makeatletter
\def\@roman#1{\romannumeral #1}
\makeatother

\newcounter{enuAlph}
\renewcommand{\theenuAlph}{\Alph{enuAlph}}

\numberwithin{equation}{section}
\renewcommand{\theequation}{\thesection.\arabic{equation}}


\theoremstyle{plain}
  \newtheorem{theorem}[equation]{Theorem}
  \newtheorem{corollary}[equation]{Corollary}
  \newtheorem{lemma}[equation]{Lemma}
  \newtheorem{mainlemma}[equation]{Main Lemma}
  \newtheorem*{mainthm}{Main Theorem}
  \newtheorem{prop}[equation]{Proposition}
  \newtheorem{clm}[equation]{Claim}
  \newtheorem{subclm}[equation]{Subclaim}
  \newtheorem{fact}[equation]{Fact}
  \newtheorem{exer}[equation]{Exercise}
  \newtheorem{question}[equation]{Question}
  \newtheorem{problem}[equation]{Problem}
  \newtheorem{conjecture}[equation]{Conjecture}
  \newtheorem{assumption}[equation]{Assumption}
    \newtheorem{hopethm}[equation]{Hopeful Theorem}
    \newtheorem{challenging}[enuAlph]{Main challenging}
    \newtheorem{hopele}[equation]{Hopeful Lemma}
    \newtheorem{discussion}[equation]{Discussion}
  \newtheorem*{thm}{Theorem}
  \newtheorem{teorema}[enuAlph]{Theorem}
  \newtheorem*{corolario}{Corollary}
\theoremstyle{definition}
  \newtheorem{definition}[equation]{Definition}
  \newtheorem{example}[equation]{Example}
  \newtheorem{remark}[equation]{Remark}
  \newtheorem{notation}[equation]{Notation}
  \newtheorem{context}[equation]{Context}

  \newtheorem*{defi}{Definition}
  \newtheorem*{acknowledgements}{Acknowledgements}

\def\sectionautorefname{Section}
\def\subsectionautorefname{Subsection}


\begin{abstract}
Let $\mathfrak{e}^\mathsf{const}_2$ be the \emph{constant evasion number}, that is, the size of the least family $F\subseteq\cantor$ of reals such that for each predictor $\pi\colon {}^{<\omega}2\to 2$ there is $x\in F$ which is not constantly predicted by $\pi$; and let $\vfa_2^\const$ be the \emph{constant prediction number}, that is, the size of the least family $\Pi_2$ of
functions $\pi\colon {}^{<\omega}2\to 2$ such that for each $x\in\cantor$ there is $\pi\in\Pi_2$ that predicts constantly $x$. In this work, we show that the constant evasion number $\mathfrak{e}_2^{\mathrm{cons}}$ and the constant prediction number $\mathfrak{v}_2^\mathsf{const}$ can be added to Cicho\'n's maximum with distinct
values.
\end{abstract}

\maketitle

\section{Introduction}\label{s0}

In this paper, we refer to a function $\pi\colon {}^{<\omega}2\to2$ as a~\textit{predictor}. Denote by $\Pi_2$ the set of predictors $\pi\colon {}^{<\omega}2\to2$. Given $\pi\in\Pi_2$ and $f\in\cantor$, say that \emph{$\pi$ predicts constantly $f$} denoted by $f\sqsubset^\pr\pi$ iff
$\exists k\in\omega$ $\forall^\infty i$ $\exists j\in[i,i+k)\colon f(j)=\pi(f{\upharpoonright}j)$.
We define the cardinal characteristics associated with $\sqsubset^\pr$.
\begin{align*}
\efrak_2^\const
&:=\min\set{|F|}{F\subseteq\cantor\setand\neg\exists \pi\in\Pi_2\ \forall f\in F\colon f\sqsubset^\pr\pi},\\
\vfa_2^\const
&:=\min\set{|S|}{S\subseteq\Pi_2\setand\forall x\in\cantor\ \exists \pi\in S\colon f\sqsubset^\pr\pi}.
\end{align*}
These cardinals are called the \emph{constant prediction number} and the \emph{constant evasion number}, respectively.

Blass~\cite{blasse} first presented the general concept of evasion and prediction concerning the Specker phenomenon in Abelian group theory. Further details on various forms of evasion can be found in~\cite[Sec.~4]{blass}. The concept of constant prediction was first attributed to Kamo, as referenced in works~\cite{kamoeva,BreIII}, and has since been explored in greater depth in publications such as~\cite{BreIII,BreShevaIV,Kadagen,BrGA}. The cardinal notation used is credited to Kada~\cite{Kaunp}.

These cardinals have connections to other cardinals. For example, Kamo~\cite{kamoeva} observed that $\vfa_2^\const\geq\cov(\Ewf)$ and $\efrak_2^\const\leq\non(\Ewf)$ where $\non(\Ewf)$ is as usual the \emph{uniformity number} of the ideal $\Ewf$ generated by the $F_\sigma$ measure zero sets on $\cantor$. He also proved that $\vfa_2^\const$ may be larger than all cardinals in Cicho\'n's diagram, and smaller than \emph{dominating number} $\dfrak$.

On the other hand, Brendle~\cite{BreIII} pointed out that $\bfrak\leq\vfa^\const_2$ and that it is relatively consistent that \emph{bounding number} $\bfrak$ is smaller than $\efrak_2^\const$.
 In contrast to the latter, it was conjectured that
\begin{conjecture}[Kada, \cite{Kaunp}]\label{a}
$\efrak^\const_2\le\dfrak$.
\end{conjecture}
Afterwards, Brendle and the third author~\cite{BreShevaIV} proved the consistency of $\efrak_2^\const<\add(\Nwf)$, where $\add(\Nwf)$ denotes \emph{additivity of the null ideal}. Dually, it was established that consistently, $\vfrak_2^\const>\cof(\Nwf)$ where \emph{$\cof(\Nwf)$ denotes the cofinality of the null ideal.}

Before stating our main theorems, we review some notation: Given a formula $\phi$,
\begin{itemize}
\item[{}] $\forall^\infty\, n<\omega\colon \phi$ means that all but finitely many natural numbers satisfy $\phi$; and
\item[{}] $\exists^\infty\, n<\omega\colon \phi$ means that infinitely many natural numbers satisfy $\phi$.
\end{itemize}
For $s\in{}^{<\omega}2=\bigcup_{n\in\omega}{}^n2$, we write $\lh (s)=\dom(s)$. If $r\in{}^{\leq\omega}2$ and $s\in{}^{\leq\omega}2$, we write $r\unlhd s$ if $r=s{\upharpoonright}\lh(r)$. Let $r\lhd s$ denote $r\unlhd s$ and $r\neq s$. A subset
$T\subseteq{}^{<\omega}2$ is called a \emph{tree} if it is closed downward, i.e., if for all $s\in T$ for all $r\unlhd s$, we have $r \in T$. 
For a~set $T\subseteq{}^{<\omega}2$, $s\in{}^{<\omega}2$,
and $m\in\omega$ denote
\begin{align*}
&[T]=\set{x\in\cantor}{\forall n\bsp x{\restriction}n\in T},\\
&T_s=\set{t\in T}{s\subseteq t\setsth{or}t\subseteq s},\\
&T^{[m]}=\set{s\in{}^{<\omega}2}
{\exists t\in T\bsp|t|=|s|\setand\forall i\in|s|\smallsetminus m\bsp
s(i)=t(i)}.
\end{align*}
Denote by $\Nwf$ and $\Mwf$ the $\sigma$-ideals of Lebesgue null sets and of meager sets in~$\cantor$, respectively; and let $\Ewf$ be the $\sigma$-ideal generated by the closed measure zero subsets of $\cantor$. It is well known that $\Ewf\subseteq\Nwf\cap\Mwf$. Even more, it was proved that $\Ewf$ is a proper subideal of $\Nwf\cap\Mwf$ (see~\cite[Lemma 2.6.1]{BJ}).
Let $\Iwf$ be an ideal of subsets of $X$ such that $\{x\}\in \Iwf$ for all $x\in X$. Throughout this paper, we demand that all ideals satisfy this latter requirement. We recall the following four \emph{cardinal invariants associated with $\Iwf$}:
\begin{align*}
\add(\Iwf)&=\min\largeset{|\Jwf|}{\Jwf\subseteq\Iwf,\ \bigcup\Jwf\notin\Iwf},\\
\cov(\Iwf)&=\min\largeset{|\Jwf|}{\Jwf\subseteq\Iwf,\ \bigcup\Jwf=X},\\
\non(\Iwf)&=\min\set{|A|}{A\subseteq X,\ A\notin\Iwf},\quad\text{and}\\
\cof(\Iwf)&=\min\set{|\Jwf|}{\Jwf\subseteq\Iwf,\ \forall A\in\Iwf\ \exists B\in \Jwf\colon A\subseteq B}.
\end{align*}
These cardinals are referred to as the \emph{additivity, covering, uniformity} and \emph{cofinality of $\Iwf$}, respectively. For $f,g\in\baire$ define
\[f\leq^*g\text{ iff } \forall^\infty n\in\omega\colon f(n)\leq g(n).\]
We recall
\begin{align*}
\bfrak &:=\min\set{|F|}{F\subseteq\baire\setand\forall g\in\baire\ \exists f\in F:f\not\leq^* g},\\
\dfrak &:=\min\set{|D|}{D\subseteq\baire\setand\forall g\in\baire\ \exists f\in D:g\leq^* f},
\end{align*}
denote the \textit{bounding number} and the \textit{dominating number}, respectively. And denote $\cfrak=2^{\aleph_0}$.

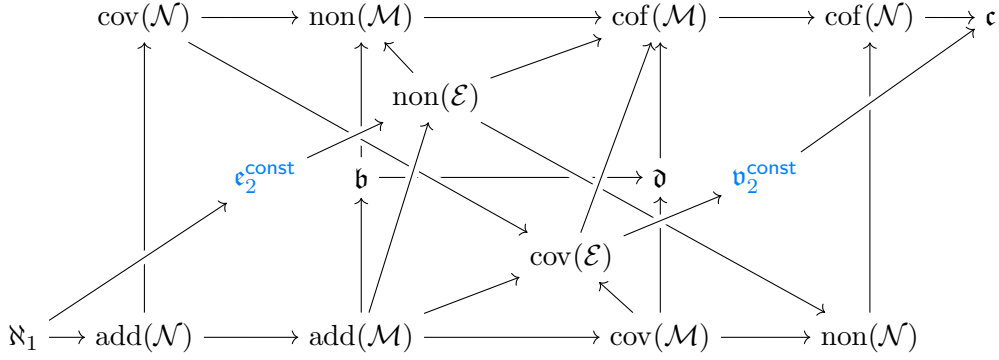
\begin{figure}[ht!]
\centering
\begin{tikzpicture}[scale=1.06]
\small{
\node (aleph1) at (-1,3) {$\aleph_1$};
\node (addn) at (0.5,3){$\add(\Nwf)$};
\node (covn) at (0.5,7){$\cov(\Nwf)$};
\node (nonn) at (9.5,3) {$\non(\Nwf)$} ;
\node (cfn) at (9.5,7) {$\cof(\Nwf)$} ;
\node (addm) at (3.19,3) {$\add(\Mwf)$} ;
\node (covm) at (6.9,3) {$\cov(\Mwf)$} ;
\node (nonm) at (3.19,7) {$\non(\Mwf)$} ;
\node (cfm) at (6.9,7) {$\cof(\Mwf)$} ;
\node (b) at (3.19,5) {$\bfrak$};
\node (d) at (6.9,5) {$\dfrak$};
\node (c) at (11,7) {$\cfrak$};
\node (e) at (2,5) {\dodger{$\efrak_2^\mathsf{const}$}};
\node (deq) at (8.2,5) {\dodger{$\mathfrak{v}_2^\mathsf{const}$}};
\node (none) at (4.12,6) {$\non(\Ewf)$};
\node (cove) at (5.8,4) {$\cov(\Ewf)$};
\draw (aleph1) edge[->] (addn)
      (addn) edge[->] (covn)
      (covn) edge [->] (nonm)
      (nonm)edge [->] (cfm)
      (cfm)edge [->] (cfn)
      (cfn) edge[->] (c);

\draw
   (addn) edge [->]  (addm)
   (addm) edge [->]  (covm)
   (covm) edge [->]  (nonn)
   (nonn) edge [->]  (cfn);
\draw (addm) edge [->] (b)
      (b)  edge [->] (nonm);
\draw (covm) edge [->] (d)
      (d)  edge[->] (cfm);
\draw (b) edge [->] (d);

\draw   (none) edge [->] (nonm)
        (none) edge [->] (cfm)
        (addm) edge [->] (cove);
      
\draw (none) edge [line width=.15cm,white,-] (nonn)
      (none) edge [->] (nonn);
      
\draw (cove) edge [line width=.15cm,white,-] (covn)
      (cove) edge [<-] (covm)
      (cove) edge [<-] (covn);

\draw (addm) edge [line width=.15cm,white,-] (none)
      (addm) edge [->] (none); 

\draw (cove) edge [line width=.15cm,white,-] (cfm)
      (cove) edge [->] (cfm);

\draw (e) edge [line width=.15cm,white,-] (none)
      (e) edge [->] (none);   

\draw (e) edge [line width=.15cm,white,-] (aleph1)
      (e) edge [<-] (aleph1);  


\draw (deq) edge [line width=.15cm,white,-] (cove)
      (deq) edge [<-] (cove);  

\draw (deq) edge [line width=.15cm,white,-] (c)
      (deq) edge [->] (c);  

}
\end{tikzpicture}
\caption{Including $\efrak^{\mathsf{const}}_2$ and $\mathfrak{v}^{\mathsf{const}}_2$ to Cicho\'n's diagram.}
\label{cichonplus}
\end{figure}

\autoref{cichonplus} illustrates the provable inequalities among the constant prediction number and the constant evasion
number, uniformity number, and covering number of $\Ewf$, and the cardinals in Cicho\'n's diagram.

Goldstern, Kellner, Mej\'ia, and the third author~\cite{GKMSe} used finitely additive measures (FAMs) along FS iterations to construct a~poset to force that
\begin{multline*}\label{cicmaxe}
 \aleph_1<\add(\Nwf)<\cov(\Nwf)<\bfrak<\efrak<\non(\Mwf)<\cov(\Mwf)<\dfrak<\\
 <\non(\Nwf)<\cof(\Nwf)<\cfrak.
 \tag{$\varotimes$}
\end{multline*}
Using the analog of this technique for closed uf-limits, Yamazoe~\cite{Ye} constructed a~poset to force~\eqref{cicmaxe}. Recently, in~\cite{Ye25} he used uf-limits on intervals (introduced by Mej\'ia~\cite{M24Anatomy}) along FS iterations to construct a~poset to force
\begin{multline*}\label{cicmaxnone}
 \aleph_1<\add(\Nwf)<\cov(\Nwf)<\bfrak<\efrak<\non(\Ewf)<\non(\Mwf)<\\
 <\cov(\Mwf)<\dfrak<\non(\Nwf)<\cof(\Nwf)<\cfrak.
 \tag{$\boxtimes$}
\end{multline*}
These two results state both possible ways to separate the cardinals in Cichoń’s diagram alongside one or more cardinal characteristics. The main challenge in both results is to separate $\efrak$ and $\non(\Ewf)$ from the rest of the cardinals in Cichoń’s diagram.

Inspired by~\eqref{cicmaxe} and \eqref{cicmaxnone}, we ask:

\begin{question}

Can it force Cichoń maximum together with\/ $\efrak_2^\const$ and\/ $\mathfrak{v}_2^\const$?
\end{question}

In this work, the prior question is answered positively. Our main results are:

\begin{teorema}\label{Thm:a0}
Let $\lambda = \lambda^{\aleph_0}$ be a cardinal and, for\/ $1\leq i\leq 5$, let $\lambda^\bfrak_i$ and $\lambda^\dfrak_i$ be uncountable regular cardinals
such that $\lambda^\bfrak_i\leq \lambda^\bfrak_j\leq \lambda^\dfrak_j\leq \lambda^\dfrak_i\leq\lambda$ for any $i\leq j$, and assume that\/ $\cof\left(([\lambda^\dfrak_1]^{<\lambda^\bfrak_1})^{\lambda^\dfrak_4}\right) = \lambda^\dfrak_1$.
Then we can construct a ccc poset forcing~\autoref{Fig:Thm:a0}.
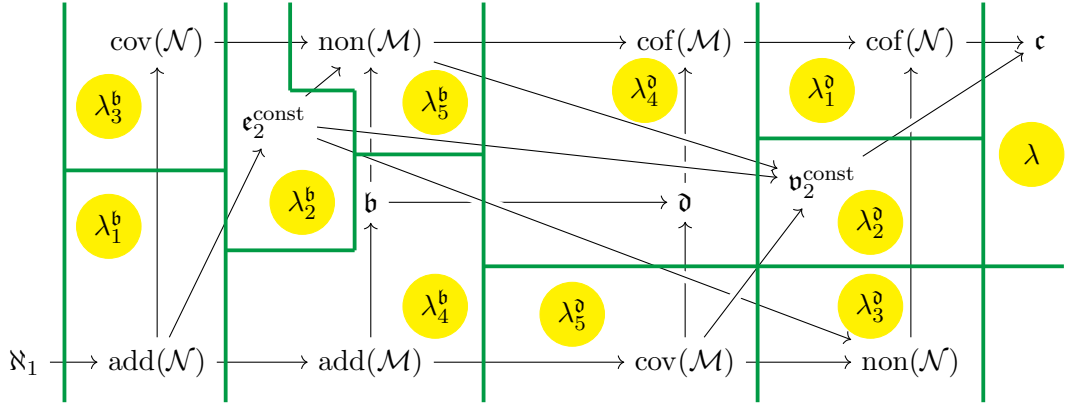
\begin{figure}[H]
\centering
\begin{tikzpicture}[scale=1.06]
\small{
\node (aleph1) at (-1.3,3) {$\aleph_1$};
\node (addn) at (0.35,3){$\add(\Nwf)$};
\node (covn) at (0.35,7){$\cov(\Nwf)$};
\node (nonn) at (9.7,3) {$\non(\Nwf)$} ;
\node (cfn) at (9.7,7) {$\cof(\Nwf)$} ;
\node (addm) at (3,3) {$\add(\Mwf)$} ;
\node (covm) at (6.9,3) {$\cov(\Mwf)$} ;
\node (nonm) at (3,7) {$\non(\Mwf)$} ;
\node (cfm) at (6.9,7) {$\cof(\Mwf)$} ;
\node (b) at (3,5) {$\bfrak$};
\node (d) at (6.9,5) {$\dfrak$};
\node (c) at (11.3,7) {$\cfrak$};
\node (cove) at (1.8,6) {$\efrak_2^\const$};
\node (vfc) at (8.6,5.25) {$\vfrak_2^\const$};
\draw (aleph1) edge[->] (addn)
      (addn) edge[->] (covn)
      (covn) edge [->] (nonm)
      (nonm)edge [->] (cfm)
      (cfm)edge [->] (cfn)
      (cfn) edge[->] (c);

\draw
   (addn) edge [->]  (addm)
   (addm) edge [->]  (covm)
   (covm) edge [->]  (nonn)
   (cove) edge [->]  (nonm)
   (nonn) edge [->]  (cfn);
\draw (addm) edge [->] (b)
      (b)  edge [->] (nonm);
\draw (covm) edge [->] (d)
      (d)  edge[->] (cfm);
\draw (b) edge [->] (d);

\draw   
        (addn) edge [->] (cove);


\draw (cove) edge [line width=.15cm,white,-] (nonn)
      (cove) edge [->] (nonn);

\draw (nonm) edge [line width=.15cm,white,-] (vfc)
      (nonm) edge [->] (vfc);

      \draw (covm) edge [line width=.15cm,white,-] (vfc)
      (covm) edge [->] (vfc);

\draw (cove) edge [line width=.15cm,white,-] (vfc)
      (cove) edge [->] (vfc);

\draw (c) edge [line width=.15cm,white,-] (vfc)
      (c) edge [<-] (vfc);

\draw[color=sug,line width=.05cm] (-0.8,5.4)--(1.2,5.4);
\draw[color=sug,line width=.05cm] (-0.8,2.5)--(-0.8,7.5);
\draw[color=sug,line width=.05cm] (1.2,2.5)--(1.2,7.5);
\draw[color=sug,line width=.05cm] (2,6.4)--(2,7.5);
\draw[color=sug,line width=.05cm] (2,6.4)--(2.8,6.4);
\draw[color=sug,line width=.05cm] (2.8,4.4)--(2.8,6.4);
\draw[color=sug,line width=.05cm] (1.2,4.4)--(2.8,4.4);
\draw[color=sug,line width=.05cm] (2.8,5.6)--(4.4,5.6);
\draw[color=sug,line width=.05cm] (4.4,2.5)--(4.4,7.5);
\draw[color=sug,line width=.05cm] (4.4,4.2)--(11.6,4.2);
\draw[color=sug,line width=.05cm] (7.8,7.5)--(7.8,2.5);
\draw[color=sug,line width=.05cm] (7.8,5.8)--(10.6,5.8);
\draw[color=sug,line width=.05cm] (10.6,7.5)--(10.6,2.5);

\draw[circle, fill=yellow,color=yellow] (-0.25,4.7) circle (0.4);
\draw[circle, fill=yellow,color=yellow] (-0.25,6.2) circle (0.4);
\draw[circle, fill=yellow,color=yellow] (3.8,3.7) circle (0.4);
\draw[circle, fill=yellow,color=yellow] (2.15,5) circle (0.4);
\draw[circle, fill=yellow,color=yellow] (3.8,6.25) circle (0.4);
\draw[circle, fill=yellow,color=yellow] (6.4,6.4) circle (0.4);
\draw[circle, fill=yellow,color=yellow] (5.5,3.6) circle (0.4);
\draw[circle, fill=yellow,color=yellow] (9.2,3.7) circle (0.4);
\draw[circle, fill=yellow,color=yellow] (8.6,6.4) circle (0.4);
\draw[circle, fill=yellow,color=yellow] (11.2,5.6) circle (0.4);
\draw[circle, fill=yellow,color=yellow] (9.2,4.75) circle (0.4);
\node at (-0.25,4.7) {$\lambda_1^\bfrak$};
\node at (-0.25,6.2) {$\lambda_3^\bfrak$};
\node at (3.8,3.7) {$\lambda_4^\bfrak$};
\node at (5.5,3.6) {$\lambda_5^\dfrak$};
\node at (2.15,5) {$\lambda_2^\bfrak$};
\node at (3.8,6.25) {$\lambda_5^\bfrak$};
\node at (6.4,6.4) {$\lambda_4^\dfrak$};
\node at (9.2,3.7) {$\lambda_3^\dfrak$};
\node at (8.6,6.4) {$\lambda_{1}^\dfrak$};
\node at (11.2,5.6) {$\lambda$};
\node at (9.2,4.75) {$\lambda_{2}^\dfrak$};
}
\end{tikzpicture}
\caption{Constellation forced in~\autoref{Thm:a0}}
\label{Fig:Thm:a0}
\end{figure}
\end{teorema}

\begin{teorema}\label{Thm:a1}
Let $\lambda = \lambda^{\aleph_0}$ be a cardinal and, for\/ $1\leq i\leq 5$, let $\lambda^\bfrak_i$ and $\lambda^\dfrak_i$ be uncountable regular cardinals
such that $\lambda^\bfrak_i\leq \lambda^\bfrak_j\leq \lambda^\dfrak_j\leq \lambda^\dfrak_i\leq\lambda$ for any $i\leq j$, and assume that\/ $\cof\left(([\lambda^\dfrak_1]^{<\lambda^\bfrak_1})^{\lambda^\dfrak_4}\right) = \lambda^\dfrak_1$.
Then we can construct a ccc poset forcing~\autoref{Fig:Thm:a1}.
\begin{figure}[H]
\centering
\begin{tikzpicture}[scale=1.06]
\small{
\node (aleph1) at (-1.3,3) {$\aleph_1$};
\node (addn) at (0.35,3){$\add(\Nwf)$};
\node (covn) at (0.35,7){$\cov(\Nwf)$};
\node (nonn) at (9.7,3) {$\non(\Nwf)$} ;
\node (cfn) at (9.7,7) {$\cof(\Nwf)$} ;
\node (addm) at (3,3) {$\add(\Mwf)$} ;
\node (covm) at (6.9,3) {$\cov(\Mwf)$} ;
\node (nonm) at (3,7) {$\non(\Mwf)$} ;
\node (cfm) at (6.9,7) {$\cof(\Mwf)$} ;
\node (b) at (3,5) {$\bfrak$};
\node (d) at (6.9,5) {$\dfrak$};
\node (c) at (11.3,7) {$\cfrak$};
\node (cove) at (1.8,6) {$\efrak_2^\const$};
\node (vfc) at (8.6,5.25) {$\vfrak_2^\const$};
\draw (aleph1) edge[->] (addn)
      (addn) edge[->] (covn)
      (covn) edge [->] (nonm)
      (nonm)edge [->] (cfm)
      (cfm)edge [->] (cfn)
      (cfn) edge[->] (c);

\draw
   (addn) edge [->]  (addm)
   (addm) edge [->]  (covm)
   (covm) edge [->]  (nonn)
   (cove) edge [->]  (nonm)
   (nonn) edge [->]  (cfn);
\draw (addm) edge [->] (b)
      (b)  edge [->] (nonm);
\draw (covm) edge [->] (d)
      (d)  edge[->] (cfm);
\draw (b) edge [->] (d);

\draw   
        (addn) edge [->] (cove);


\draw (cove) edge [line width=.15cm,white,-] (nonn)
      (cove) edge [->] (nonn);

\draw (nonm) edge [line width=.15cm,white,-] (vfc)
      (nonm) edge [->] (vfc);

      \draw (covm) edge [line width=.15cm,white,-] (vfc)
      (covm) edge [->] (vfc);

\draw (cove) edge [line width=.15cm,white,-] (vfc)
      (cove) edge [->] (vfc);

\draw (c) edge [line width=.15cm,white,-] (vfc)
      (c) edge [<-] (vfc);

\draw[color=sug,line width=.05cm] (-0.8,5.4)--(1.2,5.4);
\draw[color=sug,line width=.05cm] (-0.8,2.5)--(-0.8,7.5);
\draw[color=sug,line width=.05cm] (1.2,2.5)--(1.2,7.5);
\draw[color=sug,line width=.05cm] (2,6.4)--(2,7.5);
\draw[color=sug,line width=.05cm] (2,6.4)--(2.8,6.4);
\draw[color=sug,line width=.05cm] (2.8,4.4)--(2.8,6.4);
\draw[color=sug,line width=.05cm] (1.2,4.4)--(2.8,4.4);
\draw[color=sug,line width=.05cm] (2.8,5.6)--(4.4,5.6);
\draw[color=sug,line width=.05cm] (4.4,2.5)--(4.4,7.5);
\draw[color=sug,line width=.05cm] (4.4,4.2)--(11.6,4.2);
\draw[color=sug,line width=.05cm] (7.8,7.5)--(7.8,2.5);
\draw[color=sug,line width=.05cm] (7.8,5.8)--(10.6,5.8);
\draw[color=sug,line width=.05cm] (10.6,7.5)--(10.6,2.5);

\draw[circle, fill=yellow,color=yellow] (-0.25,4.7) circle (0.4);
\draw[circle, fill=yellow,color=yellow] (-0.25,6.2) circle (0.4);
\draw[circle, fill=yellow,color=yellow] (3.8,3.7) circle (0.4);
\draw[circle, fill=yellow,color=yellow] (2.15,5) circle (0.4);
\draw[circle, fill=yellow,color=yellow] (3.8,6.25) circle (0.4);
\draw[circle, fill=yellow,color=yellow] (6.4,6.4) circle (0.4);
\draw[circle, fill=yellow,color=yellow] (5.5,3.6) circle (0.4);
\draw[circle, fill=yellow,color=yellow] (9.2,3.7) circle (0.4);
\draw[circle, fill=yellow,color=yellow] (8.6,6.4) circle (0.4);
\draw[circle, fill=yellow,color=yellow] (11.2,5.6) circle (0.4);
\draw[circle, fill=yellow,color=yellow] (9.2,4.75) circle (0.4);
\node at (-0.25,4.7) {$\lambda_2^\bfrak$};
\node at (-0.25,6.2) {$\lambda_3^\bfrak$};
\node at (3.8,3.7) {$\lambda_4^\bfrak$};
\node at (5.5,3.6) {$\lambda_5^\dfrak$};
\node at (2.15,5) {$\lambda_1^\bfrak$};
\node at (3.8,6.25) {$\lambda_5^\bfrak$};
\node at (6.4,6.4) {$\lambda_4^\dfrak$};
\node at (9.2,3.7) {$\lambda_3^\dfrak$};
\node at (8.6,6.4) {$\lambda_{1}^\dfrak$};
\node at (11.2,5.6) {$\lambda$};
\node at (9.2,4.75) {$\lambda_{2}^\dfrak$};
}
\end{tikzpicture}
\caption{Constellation forced in~\autoref{Thm:a1}}
\label{Fig:Thm:a1}
\end{figure}
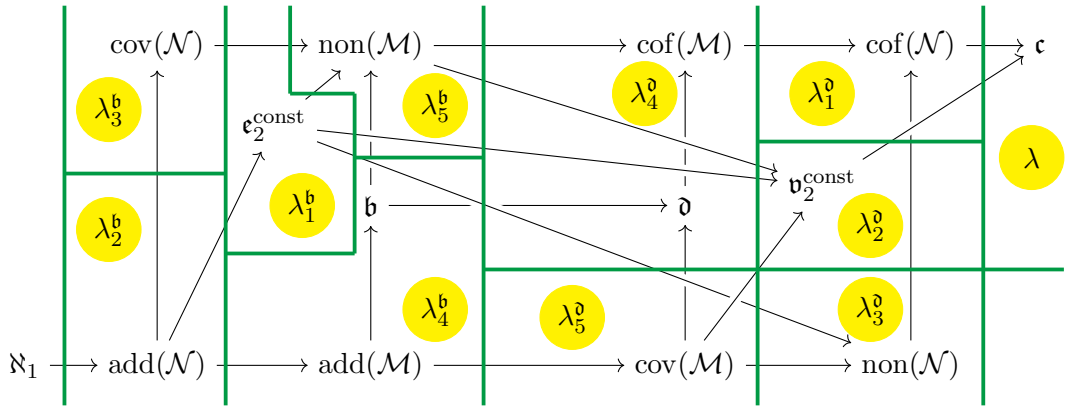
\end{teorema}

\begin{teorema}\label{Thm:a2}
Let $\lambda = \lambda^{\aleph_0}$ be a cardinal and, for\/ $1\leq i\leq 5$, let $\lambda^\bfrak_i$ and $\lambda^\dfrak_i$ be uncountable regular cardinals
such that $\lambda^\bfrak_i\leq \lambda^\bfrak_j\leq \lambda^\dfrak_j\leq \lambda^\dfrak_i\leq\lambda$ for any $i\leq j$, and assume that\/ $\cof\left(([\lambda^\dfrak_1]^{<\lambda^\bfrak_1})^{\lambda^\dfrak_4}\right) = \lambda^\dfrak_1$.
Then we can construct a ccc poset forcing~\autoref{Fig:Thm:a2}:
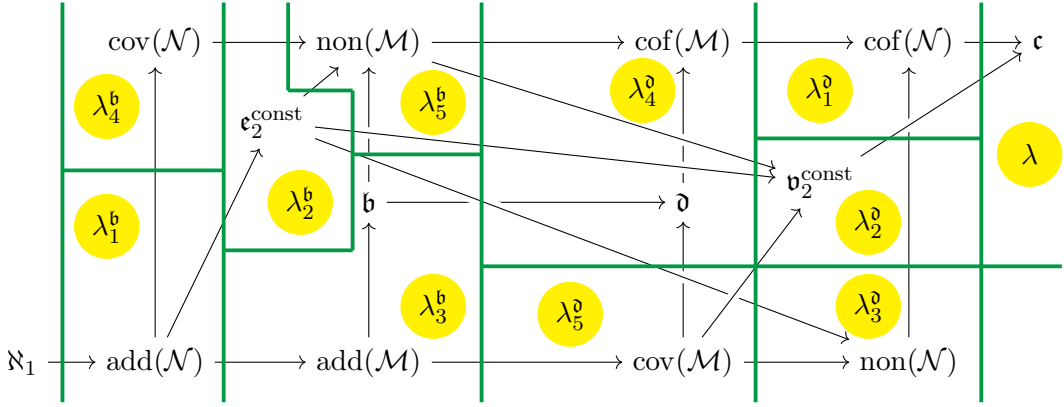
\begin{figure}[H]
\centering
\begin{tikzpicture}[scale=1.06]
\small{
\node (aleph1) at (-1.3,3) {$\aleph_1$};
\node (addn) at (0.35,3){$\add(\Nwf)$};
\node (covn) at (0.35,7){$\cov(\Nwf)$};
\node (nonn) at (9.7,3) {$\non(\Nwf)$} ;
\node (cfn) at (9.7,7) {$\cof(\Nwf)$} ;
\node (addm) at (3,3) {$\add(\Mwf)$} ;
\node (covm) at (6.9,3) {$\cov(\Mwf)$} ;
\node (nonm) at (3,7) {$\non(\Mwf)$} ;
\node (cfm) at (6.9,7) {$\cof(\Mwf)$} ;
\node (b) at (3,5) {$\bfrak$};
\node (d) at (6.9,5) {$\dfrak$};
\node (c) at (11.3,7) {$\cfrak$};
\node (cove) at (1.8,6) {$\efrak_2^\const$};
\node (vfc) at (8.6,5.25) {$\vfrak_2^\const$};
\draw (aleph1) edge[->] (addn)
      (addn) edge[->] (covn)
      (covn) edge [->] (nonm)
      (nonm)edge [->] (cfm)
      (cfm)edge [->] (cfn)
      (cfn) edge[->] (c);

\draw
   (addn) edge [->]  (addm)
   (addm) edge [->]  (covm)
   (covm) edge [->]  (nonn)
   (cove) edge [->]  (nonm)
   (nonn) edge [->]  (cfn);
\draw (addm) edge [->] (b)
      (b)  edge [->] (nonm);
\draw (covm) edge [->] (d)
      (d)  edge[->] (cfm);
\draw (b) edge [->] (d);

\draw   
        (addn) edge [->] (cove);


\draw (cove) edge [line width=.15cm,white,-] (nonn)
      (cove) edge [->] (nonn);

\draw (nonm) edge [line width=.15cm,white,-] (vfc)
      (nonm) edge [->] (vfc);

      \draw (covm) edge [line width=.15cm,white,-] (vfc)
      (covm) edge [->] (vfc);

\draw (cove) edge [line width=.15cm,white,-] (vfc)
      (cove) edge [->] (vfc);

\draw (c) edge [line width=.15cm,white,-] (vfc)
      (c) edge [<-] (vfc);

\draw[color=sug,line width=.05cm] (-0.8,5.4)--(1.2,5.4);
\draw[color=sug,line width=.05cm] (-0.8,2.5)--(-0.8,7.5);
\draw[color=sug,line width=.05cm] (1.2,2.5)--(1.2,7.5);
\draw[color=sug,line width=.05cm] (2,6.4)--(2,7.5);
\draw[color=sug,line width=.05cm] (2,6.4)--(2.8,6.4);
\draw[color=sug,line width=.05cm] (2.8,4.4)--(2.8,6.4);
\draw[color=sug,line width=.05cm] (1.2,4.4)--(2.8,4.4);
\draw[color=sug,line width=.05cm] (2.8,5.6)--(4.4,5.6);
\draw[color=sug,line width=.05cm] (4.4,2.5)--(4.4,7.5);
\draw[color=sug,line width=.05cm] (4.4,4.2)--(11.6,4.2);
\draw[color=sug,line width=.05cm] (7.8,7.5)--(7.8,2.5);
\draw[color=sug,line width=.05cm] (7.8,5.8)--(10.6,5.8);
\draw[color=sug,line width=.05cm] (10.6,7.5)--(10.6,2.5);

\draw[circle, fill=yellow,color=yellow] (-0.25,4.7) circle (0.4);
\draw[circle, fill=yellow,color=yellow] (-0.25,6.2) circle (0.4);
\draw[circle, fill=yellow,color=yellow] (3.8,3.7) circle (0.4);
\draw[circle, fill=yellow,color=yellow] (2.15,5) circle (0.4);
\draw[circle, fill=yellow,color=yellow] (3.8,6.25) circle (0.4);
\draw[circle, fill=yellow,color=yellow] (6.4,6.4) circle (0.4);
\draw[circle, fill=yellow,color=yellow] (5.5,3.6) circle (0.4);
\draw[circle, fill=yellow,color=yellow] (9.2,3.7) circle (0.4);
\draw[circle, fill=yellow,color=yellow] (8.6,6.4) circle (0.4);
\draw[circle, fill=yellow,color=yellow] (11.2,5.6) circle (0.4);
\draw[circle, fill=yellow,color=yellow] (9.2,4.75) circle (0.4);
\node at (-0.25,4.7) {$\lambda_1^\bfrak$};
\node at (-0.25,6.2) {$\lambda_4^\bfrak$};
\node at (3.8,3.7) {$\lambda_3^\bfrak$};
\node at (5.5,3.6) {$\lambda_5^\dfrak$};
\node at (2.15,5) {$\lambda_2^\bfrak$};
\node at (3.8,6.25) {$\lambda_5^\bfrak$};
\node at (6.4,6.4) {$\lambda_4^\dfrak$};
\node at (9.2,3.7) {$\lambda_3^\dfrak$};
\node at (8.6,6.4) {$\lambda_{1}^\dfrak$};
\node at (11.2,5.6) {$\lambda$};
\node at (9.2,4.75) {$\lambda_{2}^\dfrak$};
}
\end{tikzpicture}
\caption{Constellation forced in~\autoref{Thm:a2}}
\label{Fig:Thm:a2}
\end{figure}
\end{teorema}

\textbf{Method}. The proof of~\autoref{Thm:a0}, \ref{Thm:a1}, and~\ref{Thm:a2} has a similar flow as the proof of Cicho\'n's maximum~\cite{GKMS}, that is, it first forced on the left side of Cicho\'n's diagram: so we first force each one of the following statements 
\begin{enumerate}[label=\rm$\circledast_\mathrm{\Alph*}$]
\item\label{Seplef:1}
$\aleph_1<\add(\Nwf)=\lambda<\efrak_2^\const<\cov(\Nwf)<\bfrak<\non(\Mwf)<\cov(\Mwf)=\cfrak$.
\item\label{Seplef:2}
$\aleph_1<\efrak_2^\const<\add(\Nwf)<\cov(\Nwf)<\bfrak<\non(\Mwf)<\cov(\Mwf)=\cfrak$.
\item\label{Seplef:3}
$\aleph_1<\add(\Nwf)=\lambda<\efrak_2^\const<\bfrak<\cov(\Nwf)<\non(\Mwf)<\cov(\Mwf)=\cfrak$.
\end{enumerate}
and then apply the method of intersections with $\sigma$-closed models for the final of each forcing construction. The first step uses iterations with ultrafilter limits as in~\cite{GMS,GKScicmax,BCM,Car4E,BCM2,CMR2}, but we follow the two-dimensional version from~\cite{BCM}, to prove~\eqref{Seplef:1}--\eqref{Seplef:2}, except for~\eqref{Seplef:3}, to prove~\eqref{Seplef:3} uses iterations with fams limits as in~\cite{KST,CMU}; the second step's method is original by~\cite{GKMS} and is reviewed in~\cite{Brmodtec,CM22}, but our presentation is closer to~\cite{CM22}.

The main challenge in~\eqref{Seplef:1} and \eqref{Seplef:3} is to force $\add(\Nwf)=\lambda$ while iterating restrictions of the corresponding forcing notion to increase $\efrak_2^\const$ to small models. To tackle this, we introduce two new properties that behaves well to keep the additivity of the null small in forcing generics. Concretely, we introduce two properties that we call $\Pr_{\bar n}^2$ and $\lambda$-$\bar\rho$-cc for $\bar n=\seq{n_i}{i<\omega}\in\baire$, and a sequence of functions $\bar\rho=\seq{\rho_n}{n<\omega}\subseteq\baire$, respectively. These properties are something like a type of Knaster property. We prove that $\lambda$-$\bar\rho$-cc property${}\Rightarrow{}\Pr_{\bar n}^2$ property  and 

\begin{teorema}[\autoref{Thm:addN}]\label{Thm:a3}
Let $\lambda$ be a regular uncountable cardinal and $\bar n\in\baire$. Then, any iteration of length $\pi\geq\lambda$, where each iterand has the property\/ $\Pr_{\bar n}^2$, satisfies\/ $\Pr_{\bar n}^2$ and forces\/ $\add(\Nwf)\leq\lambda$ and\/ $|\pi|\leq\cof(\Nwf)$.
\end{teorema}

Further, we prove that these properties are related to the $\sigma$-$\bar\rho$-linkedness notion introduced by the authors~\cite{CRS}, which also keeps the additivity of the null small in forcing generics. This will be developed in~\autoref{sec:s3}.

On the other hand, the main work in~\eqref{Seplef:2} is to keep the evasion number $\efrak_2^\const$ small through the forcing iteration, since the other cardinal invariants are kept small without particular attention. To control $\efrak_2^\const$, in~\cite{BreShevaIV} Brendle and the third author proved that $\sigma$-$n$-linked do not increase $\efrak_2^\const$, we use this to obtain a goodness result. This is used to prove~\eqref{Seplef:2}. Details will be provided in~\autoref{Sec:sub}.

\textbf{This paper is structured as follows}. In~\autoref{sec:s1}, we review all the essentials related to relational systems, cardinal invariants, and preservation theory for cardinal characteristics. We review UF-limits and simple matrix iterations in~\autoref{sec:s2}. New tools are presented in~\autoref{sec:s3} to control the additivity of the null ideal small. We prove~\eqref{Seplef:1}--\eqref{Seplef:3} in~\autoref{sec:s4}. In \autoref{sec:subm} we review the method of forcing-intersection with $\sigma$-closed submodels and show~\autoref{Thm:a0},~\ref{Thm:a1} and~\ref{Thm:a2}. Lastly, discussions and open questions are presented in \autoref{sec:s5}. 

\section{Relational systems, cardinal invariants, and Preservation Theory for cardinal characteristics}\label{sec:s1}

We divide this section into two subsections. The first part is dedicated to recalling the basic notions and results regarding relational systems and Tukey connections introduced by Vojtáš~\cite{Vojtas} and later systematized by Blass~\cite{blass}; and the second part is devoted to reviewing the theory of the preservation of Judah and the third author~\cite{JS}, with improvements in~\cite{Br}, and it was generalized by the first author and Mej\'ia in~\cite[Sect.~4]{CM}. Furthermore, in this part one, we formalize Brendle's and the third author's result of preservation ~\cite{BreShevaIV} regarding the constant evasion number. 

\subsection{Relational systems and cardinal invariants}
\

We say that $\Rbf=\la X, Y, {\sqsubset}\ra$ is a~\textit{relational system} if it consists of two non-empty sets $X$ and $Y$ and a~relation~$\sqsubset$.
\begin{enumerate}[label=(\arabic*)]
\item A~set $F\subseteq X$ is \emph{$\Rbf$-bounded} if $\exists y\in Y$ $\forall x\in F\colon x \sqsubset y$.
\item A~set $D\subseteq Y$ is \emph{$\Rbf$-dominating} if $\forall x\in X$ $\exists y\in D\colon x \sqsubset y$.
\end{enumerate}

We associate two cardinal characteristics with this relational system $\Rbf$:
\begin{align*}
\bfrak(\Rbf)&:=\min\set{|F|}{\text{$F\subseteq X$ is $\Rbf$-unbounded}},
&&\text{the \emph{unbounding number of\/ $\Rbf$}, and}\\
\dfrak(\Rbf)&:=\min\set{|D|}{\text{$D\subseteq Y$ is $\Rbf$-dominating}},
&&\text{the \emph{dominating number of\/ $\Rbf$}.}
\end{align*}
The dual of $\Rbf$ is defined by $\Rbf^\perp:=\la Y, X, {\sqsubset^\perp}\ra$ where $y\sqsubset^\perp x$ iff $x\not\sqsubset y$. Note that $\bfrak(\Rbf^\perp)=\dfrak(\Rbf)$ and $\dfrak(\Rbf^\perp)=\bfrak(\Rbf)$.

The cardinal $\bfrak(\Rbf)$ may be undefined, in which case we write $\bfrak(\Rbf) = \infty$, as well as for $\dfrak(\Rbf)$. Concretely, $\bfrak(\Rbf) = \infty$ iff $\dfrak(\Rbf) =1$; and $\dfrak(\Rbf)= \infty$ iff $\bfrak(\Rbf) =1$.

\begin{definition}\label{def:defrel}
We say that $\Rbf=\la X, Y, {\sqsubset}\ra$ is a \textit{definable relational system of the reals} if both $X$ and $Y$ are non-empty and analytic in Polish spaces $Z$ and $W$, respectively, and $\sqsubset$ is analytic in $Z\times W$.\footnote{In general, we need that $X$, $Y$ and $\sqsubset$ are definable and that the statements $\lqq x\in X\rqq$, $\lqq y\in Y\rqq$ and $\lqq x\sqsubset y\rqq$ are absolute for the arguments we need to carry on.} The interpretation of $\Rbf$ in any model corresponds to the interpretation of $X$, $Y$ and $\sqsubset$.
\end{definition}

As in~\cite{GKMS,CM22}, we also look at relational systems given by directed preorders.

\begin{definition}\label{examSdir}
We say that $\la S,{\leq_S}\ra$ is a \emph{directed preorder} if it is a preorder (that is, $\leq_S$~is a reflexive and transitive relation in $S$) such that
\[\forall x, y\in S\ \exists z\in S\colon x\leq_S z\setand y\leq_S z.\]
A directed preorder $\la S,{\leq_S}\ra$ is seen as the relational system $S=\la S, S,{\leq_S}\ra$, and their associated cardinal characteristics are indicated by $\bfrak(S)$ and $\dfrak(S)$. The cardinal $\dfrak(S)$ is actually the \emph{cofinality of $S$}, typically denoted by $\cof(S)$ or $\cf(S)$.
\end{definition}

Note that $\leq^*$ is a directed preorder on $\baire$, where $x\leq^* y$ means $\forall^\infty n<\omega\colon x(n)\leq y(n)$. We think of $\baire$ as the relational system with the relation $\leq^*$. Then $\bfrak:=\bfrak(\baire)$ and $\dfrak:=\dfrak(\baire)$.

Relational systems can also characterize the cardinal characteristics associated with an ideal.

\begin{example}\label{exm:Iwf}
For $\Iwf\subseteq\pts(X)$, define the relational systems:
\begin{enumerate}[label=\rm(\arabic*)]
\item
$\Iwf:=\la\Iwf,\Iwf,{\subseteq}\ra$, which is a~directed partial order when $\Iwf$ is closed under unions (e.g.\ an ideal).

\item
$\Cv_\Iwf:=\la X,\Iwf,{\in}\ra$.
\end{enumerate}
Whenever $\Iwf$ is an ideal on $X$,
\begin{multicols}{2}
\begin{enumerate}[label=\rm(\alph*)]
\item
$\bfrak(\Iwf)=\add(\Iwf)$. 

\item
$\dfrak(\Iwf)=\cof(\Iwf)$. 

\item
$\dfrak(\Cv_\Iwf)=\cov(\Iwf)$. 

\item
$\bfrak(\Cv_\Iwf)=\non(\Iwf)$. 
\end{enumerate}
\end{multicols}
\end{example}

The Tukey connection is a useful tool for establishing relationships between cardinal characteristics.
Let $\Rbf=\la X,Y,{\sqsubset}\ra$ and $\Rbf'=\la X',Y',{\sqsubset'}\ra$ be two relational systems. We say that
$(\Psi_-,\Psi_+)\colon\Rbf\to \Rbf'$
is a~\emph{Tukey connection from $\Rbf$ into $\Rbf'$} if
$\Psi_-\colon X\to X'$ and $\Psi_+\colon Y'\to Y$ are functions such that
\[
\forall x\in X\ \forall y'\in Y'\colon
\Psi_-(x) \sqsubset' y' \Rightarrow x \sqsubset \Psi_+(y').
\]
The \emph{Tukey order} between relational systems is defined by
$\Rbf\leqT \Rbf'$ iff there is a~Tukey connection from $\Rbf$ into $\Rbf'$. \emph{Tukey equivalence} is defined by $\Rbf\eqT \Rbf'$ iff $\Rbf\leqT \Rbf'$ and $\Rbf'\leqT \Rbf$. Notice that $\Rbf\leqT \Rbf'$ implies $(\Rbf')^\perp\leqT \Rbf^\perp$, $\dfrak(\Rbf)\leq\dfrak(\Rbf')$ and $\bfrak(\Rbf')\leq\bfrak(\Rbf)$.

To close this section, note that we may also consider the constant prediction and evasion numbers for ${}^{\omega}K$. This is done in exactly the same fashion as for the Cantor $\cantor$: for this propose let $2\leq K\leq\omega$, we say that \emph{$\pi\colon {}^{<\omega}K\to K$ predicts constantly $f\in {}^{\omega}K$} denoted by $f\sqsubset^\pr \pi$ iff
$\exists k\in\omega$ $\forall^\infty i$ $\exists j\in[i,i+k)\colon f(j)=\pi(f{\upharpoonright}j)$.

We define the cardinal characteristics associated with $\sqsubset^\pr$.
\begin{align*}
\efrak_K^\const&:=\min\set{|F|}{F\subseteq {}^{\omega}K\setand\neg\exists \pi\in\Pi_K\ \forall f\in F\colon f\sqsubset^\pr\pi},\\
\vfa_K^\const&:=\min\set{|S|}{S\subseteq\Pi_K\setand\forall x\in {}^{\omega}K\ \exists \pi\in S\colon f\sqsubset^\pr\pi}.
\end{align*}
Let $\efrak^\const:=\efrak^\const_\omega$ and $\vfrak^\const:=\vfrak^\const_\omega$. Let $\Esf_K:=\la{}^{\omega}K,\Sigma_K,{\sqsubset^\pr}\ra$

\begin{fact}\label{t2}
$\Esf_2\leqT\Esf_K\leqT\Esf_\omega$. As a consequence, $\efrak^\const\leq\efrak_K^\const\leq\efrak_2^\const$ and\/ $\vfrak^\const_2\leq\vfrak_K^\const\leq\vfrak^\const$. 
\end{fact}

Regarding consistency, Brendle~\cite{BreIII} proved the consistency of $\efrak^\const_2>\efrak^\const$ and $\efrak^\const>\bfrak$; and Kamo showed $\vfa_2<\vfa^\const$~\cite{kamoeva} and $\vfa^\const<\dfrak$~\cite{kamopred}.

\subsection{Preservation properties for cardinal characteristics}\label{Sec:sub}

\begin{definition}\label{DefRelSys}
Let $\Rbf=\la X,Y,{\sqsubset}\ra$ be a relational system and let $\theta$ be a cardinal.
\begin{enumerate}[label=\rm(\arabic*)]
\item
For a set $M$,
\begin{enumerate}[label=\rm(\roman*)]
\item
An object $y\in Y$ is \textit{$\Rbf$-dominating over $M$} if $x\sqsubset y$ for all $x\in X\cap M$.
\item
An object $x\in X$ is \textit{$\Rbf$-unbounded over $M$} if it $\Rbf^\perp$-dominating over $M$, that is, $x\not\sqsubset y$ for all $y\in Y\cap M$.
\end{enumerate}
\item
A family $\set{x_i}{i\in I}\subseteq X$ is \emph{strongly $\theta$-$\Rbf$-unbounded} if
$|I|\geq\theta$ and, for any $y\in Y$, $|\set{i\in I }{x_i\sqsubset y}|<\theta$.
\end{enumerate}
\end{definition}

The existence of strongly unbounded families is equivalent to a Tukey connection.

\begin{lemma}[{\cite[Lem.~1.16]{CM22}}]\label{unbT}
Let\/ $\Rbf=\la X,Y,{\sqsubset}\ra$ be a relational system, $\theta$ be an infinite cardinal, and $I$ be a set of size\/ ${\geq}\theta$.  Then there exists a strongly $\theta$-$R$-unbounded family $\set{x_i}{i\in I}$ iff $\Cbf_{[I]^{<\theta}}\leqT R$. In particular, $\bfrak(\Rbf)\leq\non([I]^{<\theta})=\theta$ and
\[
\dfrak(\Rbf)\geq\cov([I]^{<\theta})=
\begin{cases}
|I|,&\text{if $\theta<|I|$,}\\
\cf(\theta),&\text{if $\theta=|I|$.}
\end{cases}
\]
\end{lemma}

Hence, $\cov([I]^{<\theta})=|I|$ when $\theta$ is regular. We below look at the following type of well-defined relational systems.

\begin{definition}\label{b19}
We say that $\Rbf=\langle X,Y,\sqsubset\rangle$ is a \textit{Polish relational system (Prs)} if
\begin{enumerate}[label=\rm(\arabic*)]
\item
$X$ is a Perfect Polish space,
\item
$Y$ is analytic in some Polish space $Z$, and
\item\label{b19:3}
${\sqsubset}=\bigcup_{n<\omega}{\sqsubset_n}$ where $\langle\sqsubset_n: n<\omega\rangle$ is some increasing sequence of closed subsets of $X\times Z$ such that, for any $n<\omega$ and for any $y\in Y$, $({\sqsubset_n})^{y}=\set{x\in X}{x\sqsubset_n y}$ is closed nowhere dense in $X$.
\end{enumerate}
\end{definition}

\begin{remark}\label{b18}
By~\autoref{b19}~\ref{b19:3}, $\la X,\Mwf(X),{\in}\ra\leqT \Rbf$ where $\Mwf(X)$ denotes the $\sigma$-ideal of meager subsets of $X$. Therefore, $\bfrak(\Rbf)\leq \non(\Mwf)$ and $\cov(\Mwf)\leq\dfrak(\Rbf)$.
\end{remark}

From now on, fix a Prs $\Rbf=\langle X,Y,{\sqsubset}\rangle$ and an infinite cardinal $\theta$.

\begin{definition}[Judah and Shelah {\cite{JS}}, Brendle~{\cite{Br}}]\label{b17}
A forcing notion $\Por$ is \textit{$\theta$-$\Rbf$-good} if, for any $\Por$-name $\dot{h}$ for a member of $Y$, there is a non-empty set $H\subseteq Y$ (in the ground model) of size ${<}\theta$ such that, for any $x\in X$, if $x$ is $\Rbf$-unbounded over $H$ then $\Vdash\lqq x\not\sqsubset \dot{h}\rqq$.

We say that $\Por$ is \textit{$\Rbf$-good} if it is $\aleph_1$-$\Rbf$-good.
\end{definition}

\begin{remark}
Note that $\theta<\theta'$ implies that any $\theta$-$\Rbf$-good poset is $\theta'$-$\Rbf$-good. Also, 
``$\theta$-$\Rbf$-good'' is a hereditary forcing property.
\end{remark}

We now provide some examples of good forcing notions. In general, ``small'' forcing notions are automatically good.

\begin{lemma}[{\cite[Lemma~4.10]{CM}}]\label{b22}
If $\theta$ is a regular cardinal then any poset of size ${<}\theta$
is $\theta$-$\Rbf$-good. In particular, Cohen forcing\/ $\Cor$ is\/ $\Rbf$-good.
\end{lemma}

\begin{example}\label{b16}
The following are Prs that describe the cardinal characteristics of Cicho\'n's diagram.
\begin{enumerate}[label=(\arabic*)]
\item
Define the relational system $\Mg:=\la\cantor,\Xi,\in^{\bullet}\ra$ where
\[\Xi := \set{f\colon{}^{<\omega}2\to{}^{<\omega}2}{\forall s \in {}^{<\omega}2\colon s \subseteq f(s)}\]
and $x\in^{\bullet} f$ iff $|\set{s\in {}^{<\omega}2}{x \supseteq f(s)}|<\aleph_0$. This is a Prs and $\Mg\eqT\Cv_\Mwf$. Hence $\bfrak(\Mg)=\non(\Mwf)$ and $\dfrak(\Mg)=\cov(\Mwf)$.

\item
The relational system $\baire:=\la\baire,\baire,{\leq^*}\ra$ is already Polish. Typical examples of $\baire$-good sets are the standard eventually real forcing $\Eor$ and random forcing $\Bor$. More generally, $\sigma$-$\Fr$-linked posets are $\D$-good (see~\cite{mejvert,BCM}).

\item\label{b16:2}
Define $\Omega_n:=\set{a\in [{}^{<\omega}2]^{<\aleph_0}}{\Lb(\bigcup_{s\in a}[s])\leq2^{-n}}$ (endowed with the discrete topology) where $\Lb$ is the Lebesgue measure on $\cantor$. Put $\Omega:=\prod_{n<\omega}\Omega_n$ with the product topology, which is a perfect Polish space. For every $x\in \Omega$ denote $N_{x}^{*}:=\bigcap_{n<\omega}\bigcup_{s\in x(n)}[s]$, which is clearly a Borel null set in $2^{\omega}$.

Define the Prs $\Cn:=\la\Omega,\cantor,{\sqsubset}\ra$ where $x\sqsubset z$ iff $z\notin N_{x}^{*}$. Recall that any null set in $\cantor$ is a subset of $N_{x}^{*}$ for some $x\in \Omega$, so $\Cn\eqT\Cv_\Nwf^\perp$. Hence, $\bfrak(\Cn)=\cov(\Nwf)$ and $\dfrak(\Cn)=\non(\Nwf)$. Any $\mu$-centered poset is $\mu^+$-$\Cn$-good (\cite{Br}). In particular, $\sigma$-centered posets are $\Cn$-good.



\item\label{b16:3} For $\Hwf\subseteq\omega^\omega$ denote $\Lc^*_\Hwf:=\langle\omega^\omega,\Swf(\omega,\Hwf),\in^{*}\rangle$ where \[\Swf(\omega, \Hwf):=\{\varphi:\omega\to[\omega]^{<\aleph_0}:\exists{h\in\Hwf}\,\forall{i<\omega}(|\varphi(i)|\leq h(i))\}.\]
    Any $\mu$-centered poset is $\mu^+$-$\Lc^*_\Hwf$-good (see~\cite{Br,JS}) so, in particular, $\sigma$-centered posets are $\Lc^*_\Hwf$-good. 
     
     If $\Hwf$ is countable and non-empty then $\Lc^*_\Hwf$ is a Prs because $\Swf(\omega,\Hwf)$ is $F_\sigma$ in $([\omega]^{<\aleph_0})^\omega$. In addition, if $\Hwf$ contains a function that goes to infinity, then $\Lc^*_\Hwf\eqT\Nwf$, so
$\bfrak(\Lc^*_\Hwf)=\add(\Nwf)$ and $\dfrak(\Lc^*_\Hwf)=\cof(\Nwf)$ see~\cite[Thm.~2.3.9]{BJ} (or see also~\cite[Thm.~4.2]{CMlocalc}). Moreover, if for any $h\in\Hwf$ there is an $h^\prime\in\Hwf$ such that $\frac{h}{h^\prime}$ coverge to $0$, then any Boolean algebra with a strictly
positive finitely additive measure is $\Lc^*_\Hwf$-good (\cite{Ka}). In particular, any
subalgebra of random forcing is $\Lc^*_\Hwf$-good.

\item\label{b16:5}
Let $\Mbf := \la 2^\omega,\Ior\times 2^\omega,{\sqsubm}\ra$ where
\[x \sqsubm (I,y)\text{ iff }\forall^\infty n\colon x\frestr I_n \neq y\frestr I_n.\]
This is a Polish relational system and $\Mbf\eqT\Cv_\Mwf$ (by Talagrand~\cite{Tal98}, see e.g.~\cite[Prop.~13]{BWS}).

Note that, whenever $M$ is a transitive model of $\thzfc$, $c\in 2^\omega$ is a~Cohen real over $M$ iff $c$ is $\Mbf$-unbounded over $M$.

\end{enumerate}
\end{example}

\begin{example}\label{b30}
Fix $k<\omega$, denote by $\Sigma_2^k$ the collection of maps of $\sigma\colon\bigcup_{i<\omega}2^{ik}\to2^{k}$. Let $\Esf_2^k=\la\cantor,\Sigma_2^k,{\sqsubset^\star}\ra$ where $x\sqsubset^\star\sigma$ iff $\forall^\infty i\colon x{\upharpoonright}[ik,(i+1)k)\neq\sigma(x{\upharpoonright}ik)$. It can be proved that the relational system $\Esf_2^k$ is a Prs.
\end{example}

Notice that

\begin{fact}\label{cxn:e_2}
Fix $k<\omega$. Then $\Esf_2^k\leqT\Esf_2:=\la\cantor,\Sigma_2,{\sqsubset^\pr}\ra$. In particular, $\efrak_2^\const=\bfrak(\Esf_2)\leq\bfrak(\Esf_2^k)$ and\/ $\dfrak(\Esf_2^k)\leq\dfrak(\Esf_2)=\vfrak_2^\const$.
\end{fact}

\begin{varproof}
Note that $\Psi_-\colon\cantor\to\cantor$, defined as the identity map, and $\Psi_+\colon\Sigma_2\to\Sigma_2^k$ is defined as $\Psi_+(\pi)(\sigma)=$ the unique $\tau\in{}^{k}2$
such that $\pi$ predicts $\sigma\char 94\tau$ incorrectly on the whole interval $[ik,(i + 1)k)$ where $|\sigma| = ik$. It is clear that
$(\Psi_-,\Psi_+)$ is the required Tukey connection.
\end{varproof}

\begin{lemma}[{\cite{BreShevaIV}}]\label{c10}
Let $k<\omega$. If\/ $\Por$ is $\sigma$-$2^k$-linked, then\/ $\Por$ is\/ $\Esf_2^k$-good.
\end{lemma}

\begin{varproof}
Suppose that $\Por$ is $\sigma$-$2^k$-linked witnessed by $\seq{P_n}{n\in\omega}$. Let $\dot\phi$
is a $\Por$–name for a function $\bigcup_{i<\omega}{}^{ik}2\to{}^{k}2$. For each $n\in\omega$ define $\psi_n\colon\bigcup_{i<\omega}{}^{ik}2\to{}^{k}2$
such that, for each $\sigma\in{}^{ik}2$
\[\text{$\psi_n(\sigma)$ is a $\tau$ such that no $p\in P_n$ forces 
$\dot\phi(\sigma)=\tau$.}\]
Such a $\tau$ clearly exists. Otherwise, for each $\tau\in {}^{k}2$ we could
find $p_\tau\in P_n$ forcing $\dot\phi(\sigma)=\tau$. Since $P_n$ is $2^k$-linked, the $p_\tau$ would have
a common extension which would force $\dot\phi(\sigma)\in {}^{k}2$, a contradiction. Let $H:=\set{\psi_n}{n\in\omega}$.

Now assume that $x\in\cantor$ such that for all $n\in\omega\colon x\not\sqsubset^\star\psi_n$ and show $\Vdash\lqq x\not\sqsubset^\star\dot\phi\rqq$.
Fix $i_0$ and $p \in \Por$. Then there is $n$ such that
$p \in P_n$. We can find $i \geq i_0$ such that $\psi_n(x{\restriction} ik) = x{\restriction}[ik,(i + 1)k)$. By
definition of $\psi_n$, there is $q \leq p$ such that $q\Vdash\lqq\dot\phi(x{\restriction} ik) = \psi_n(x{\restriction} ik)\rqq$. Thus
$q\Vdash\lqq\dot\phi(x{\restriction} ik) = x{\restriction}[ik,(i + 1)k)\rqq$, as required.
\end{varproof}

As an immediate consequence, we get

\begin{corollary}\label{c11}
Any $\sigma$-centered posets is\/ $\Esf_2^k$-good.
\end{corollary}

The following describes how good posets are preserved along FS iterations.

\begin{theorem}[{\cite[Sec.~4]{BCM2}}]\label{b4}
Let $\theta$ be an uncountable regular cardinal and \/ $\seq{\Por_\xi,\Qnm_\xi}{\xi<\pi}$ be a FS iteration such that, for $\xi<\pi$, $\Por_\xi$ forces that\/ $\Qnm_\xi$ is a non-trivial $\theta$-cc $\theta$-$\Rbf$-good poset.
Let\/ $\set{\gamma_\alpha}{\alpha<\delta}$ be an increasing enumeration of\/ $0$ and all limit ordinals smaller than $\pi$ (note that $\gamma_\alpha=\omega\alpha$), and for $\alpha<\delta$ let $\dot c_\alpha$ be a\/~$\Por_{\gamma_{\alpha+1}}$-name of a~Cohen real in $X$ over $V_{\gamma_\alpha}$.

Then\/ $\Por_\pi$ is $\theta$-$\Rbf$-good. Moreover,
if $\pi\geq\theta$ then\/ $\Cv_{[\pi]^{<\theta}}\leqT\Rbf$, $\bfrak(\Rbf)\leq\theta$ and\/ $|\pi|\leq\dfrak(\Rbf)$.
\end{theorem}

The Cohen reals added along an iteration are usually used as witnesses for Tukey connections, as they form strong witnesses. For example:

\begin{lemma}[{\cite[Lemma~4.14]{CM}}]\label{lem:strongCohen}
Let $\mu$ be a cardinal with uncountable cofinality and let\/ $\seq{\Por_{\alpha}}{\alpha<\mu}$ be a\/ $\subsetdot$-increasing sequence of\/ $\cf(\mu)$-cc posets and let\/ $\Por_\mu=\limdir_{\alpha<\mu}\Por_{\alpha}$. If\/ $\Por_{\alpha+1}$ adds a~Cohen real $\dot{c}_\alpha\in X$ over $V^{\Por_\alpha}$ for any $\alpha<\mu$, then\/ $\Por_{\mu}$ forces that\/ $\set{\dot{c}_\alpha}{\alpha<\mu}$ is a strongly $\cf(\mu)$-$\Rbf$-unbounded family. In particular, $\Por_\mu$ forces that $\mu\leqT\Cv_\Mwf\leqT \Rbf$.
\end{lemma}

Concerning FS iterations, we recall the following useful result to force statements of the form $\Rbf\leqT\Cv_{[X]^{<\theta}}$ for the following type of relational systems. Before stating recall $y\in Y$ is said \emph{$\Rbf$-dominating over a set $M$} if $\forall x \in X\cap M\colon x \sqsubset y$.

\begin{theorem}[{\cite[Thm.~2.12]{CM22}}]\label{b12}
Let\/ $\Rbf=\la X,Y,{\sqsubset}\ra$ be a definable relational system of the reals, $\theta$ an uncountable regular cardinal, and let\/ $\Por_\pi=\seq{\Por_\xi,\Qnm_\xi}{\xi<\pi}$ be a FS iteration of $\theta$-cc posets with\/ $\cf(\pi)\geq\theta$. Assume that, for all $\xi<\pi$ and any $A\in[X]^{<\theta}\cap V_\xi$, there is some $\eta\geq\xi$ such that\/ $\Qnm_\eta$ adds an\/ $\Rbf$-dominating real over $A$. Then\/ $\Por_\pi$ forces\/ $\Rbf\leqT\Cv_{[X]^{<\theta}}$, i.e.\ $\theta\leq\bfrak(\Rbf)$.
\end{theorem}

\section{Reviewing UF-limits and simple matrix iterations}\label{sec:s2}

This section consists of two parts. The former part is devoted to presenting the notion of ultrafilter limits for forcing notions from~\cite{GMS,BCM,CMR2}. These properties will be applied in the proof of our consistency results in~\autoref{sec:s4}. In addition to this, we introduce a forcing notion to increase $\efrak_K^\const$ and show that this has ultrafilter limits; in the latter part, we examine the matrix iterations with the ultrafilter method from~\cite{BCM}, which, as we mentioned, will be used to prove~\autoref{Thm:a0} and~\ref{Thm:a1}.

\subsection{UF-limits}
\

Given a~poset $\Por$, the $\Por$-name $\dot{G}$ usually denotes the canonical name of the $\Por$-generic set. If $\bar{p}=\seq {p_n}{n<\omega}$ is a sequence in $\Por$, denote by $\dot{W}_\Por(\bar{p})$ the $\Por$-name of $\set{n<\omega}{p_n\in\dot{G}}$. When the forcing is understood from the context, we just write $\dot{W}(\bar{p})$.

\begin{definition}[{\cite{GMS,BCM,CMR2}}]\label{Def:GMS}
Let $\Por$ be a~poset, $D\subseteq\pts(\omega)$ a non-principal ultrafilter, and $\mu$ an infinite cardinal.
\begin{enumerate}[label=\rm(\arabic*)]
\item\label{GMS1} A set $Q\subseteq \Por$ has \emph{$D$-limits} if there is a function $\lim^{D}\colon Q^\omega\to \Por$ and a $\Por$-name $\dot D'$ of an ultrafilter extending $D$ such that, for any $\bar q = \seq{q_i}{i<\omega}\in Q^\omega$,
\[\lim\nolimits^{D}\bar q\Vdash\lqq\dot{W}(\bar{q})\in \dot D'\,\rqq.\]

\item A set $Q\subseteq \Por$ has \emph{uf-limits} if it has $D$-limits for any ultrafilter $D$.

\item The poset $\Por$ is \emph{$\mu$-$D$-$\lim$-linked} if $\Por = \bigcup_{\alpha<\mu}Q_\alpha$ where each $Q_\alpha$ has $D$-limits. We say that $\Por$ is \emph{uniformly $\mu$-$D$-$\lim$-linked} if, additionally, the $\Por$-name $\dot D'$ from~\ref{GMS1} only depends on $D$ (and not on $Q_\alpha$, although we have different limits for each $Q_\alpha$).

\item The poset $\Por$ is \emph{$\mu$-uf-$\lim$-linked} if $\Por = \bigcup_{\alpha<\mu}Q_\alpha$ where each $Q_\alpha$ has uf-limits. We say that $\Por$ is \emph{uniformly $\mu$-uf-$\lim$-linked} if, additionally, for any ultrafilter $D$ on $\omega$, the $\Por$-name $\dot D'$ from~\ref{GMS1} only depends on $D$.
\end{enumerate}
\end{definition}

To not add dominating reals, we have the following weaker notion.

\begin{definition}[{\cite{mejvert}}]\label{Def:Fr}
Let $\Por$ be a~poset and $F$ a filter on $\omega$. A set $Q\subseteq \Por$ is \emph{$F$-linked} if, for any $\bar p=\seq{p_n}{n<\omega} \in Q^\omega$, there is some $q\in \Por$ forcing that $F\cup \{\dot{W}(\bar p)\}$ generates a~filter on $\omega$.
We say that $Q$ is \emph{uf-linked (ultrafilter-linked)} if it is $F$-linked for any filter~$F$ on $\omega$ containing \emph{Frechet filter} $\Fr:=\set{\omega\menos a}{a\in[\omega]^{<\aleph_0}}$.

For an infinite cardinal $\mu$, $\Por$ is \emph{$\mu$-$F$-linked} if $\Por = \bigcup_{\alpha<\mu}Q_\alpha$ for some $F$-linked $Q_\alpha$ ($\alpha<\mu$). When these $Q_\alpha$ are uf-linked, we say that $\Por$ is \emph{$\mu$-uf-linked}.
\end{definition}

It is clear that any uf-$\lim$-linked set $Q\subseteq \Por$ is uf-linked, which implies $\Fr$-linked.

\begin{theorem}[{\cite{mejvert}}]\label{mej:uf}
Any $\mu$-$\Fr$-linked poset is $\mu^+$-$\omega^\omega$-good.
\end{theorem}

\begin{example}\label{exm:ufl}
The following are the instances of $\mu$-uf-lim-linked posets that we use in our applications.
\begin{enumerate}[label=\rm (\arabic*)]
\item
Any poset of size $\mu$ is uniformly $\mu$-uf-lim-linked (because singletons are uf-lim-linked). In particular, Cohen
forcing is uniformly $\sigma$-uf-lim-linked.
\item
$\Eor$ is uniformly $\sigma$-uf-lim-linked (see~\cite{GMS,BCM}).
\item
$\Bor$ is $\sigma$-uf-linked~\cite[Lem.~3.29 \& Lem.~5.5]{mejvert} but may not be $\sigma$-uf-$\lim$-linked (cf.~\cite[Rem.~3.10]{BCM}).
\end{enumerate}
\end{example}

Below, we introduce a forcing notion to increase $\efrak_K^\const$. This generalizes a forcing notion of Brendle~\cite{BreIII}.

\begin{definition}\label{forc:Pr_K}
Let $2\leq K<\omega$. The conditions in $\Por^K$ are triples
$p=(n_p,\pi_p,F_p)=(n,\pi,F)$ such that
\begin{itemize}
\item
$n\in\omega$,
$\pi\colon{}^{\leq n}K\to\{0,1\}$,
$F\in[{}^{\omega}K]^{<\omega}$,
\item
$\eta{\restriction}n\neq\eta'{\restriction}n$ for any
different $\eta,\eta'\in F_p$, and
\item
$\eta(n)=\pi(\eta{\restriction}n)$ for each $\eta\in F$.
\end{itemize}
$\Por^K$ is partially ordered by
\begin{align*}
q\leq p\Leftrightarrow{}
&\text{$n_p\leq n_q$, $\pi_p\subseteq\pi_q$, $F_p\subseteq F_q$, and}\\
&\forall\eta\in F_p\ \forall\ell\in[n_p,n_q)\ \exists j\in\{\ell,\ell+1\}\ \eta(j)=\pi_q(\eta{\restriction}j).
\end{align*}
\end{definition}

Let $G\subseteq\Por^K$ be a $\Por^K$-generic over $V$.
Then, in $V[G]$, define
\[
\pi_G=\bigcup\set{\pi_p}{p\in G}.
\]

By a standard density argument, it is not hard to see that $\pi_G\in\Predictors_K$.

\begin{fact}
Every real $\eta\in{}^{\omega}K\cap V$ is predicted by $\pi_G$.
\end{fact}

We below drawn our attention to prove that $\Por^K$ has UF-limits.

\begin{lemma}\label{Pr_K-uf}
Let $K$ as in~\autoref{forc:Pr_K}. Then\/ $\Por^K$ is uniformly $\sigma$-uf-lim-linked.
\end{lemma}

\begin{varproof}
For $k$, $\pi$ as in~\autoref{forc:Pr_K}, and $\bar\phi=\{\phi_0,\phi_1,\ldots,\phi_{i-1}\}\subseteq{}^{k}K$ let
\[Q_{k,\pi,\bar\phi}=\set{(k,\pi,F)\in\Pre}{|F|=i\setand\forall f\in F\ \exists j<i\colon f{\upharpoonright}k=\phi_j}.\]
Let $D$ be a non-principal ultrafilter on $\omega$, and $\bar{p}=\seq{p_m}{m\in\omega}$ be a sequence in $Q_{k,\pi,\bar\phi}$ with $p_m = (k,\pi,F_m)$ and $F_m=\set{f_{m,j}}{j<i}$ where $f_{m,j}{\upharpoonright}k=\phi_j$. We define $\lim^D\bar p$. Considering the lexicographic order $\leq_{\mathrm{lex}}$ of ${}^{\omega}K$, and let $\set{f_{n,j}}{j<i}$ be a $\leq_{\mathrm{lex}}$-increasing enumeration of $F_m$. For each $j<i$, define $f_j:=\lim_m^D f_{m,j}$ in ${}^{\omega}K$ where
\[\text{$f_j(n)=\ell$ iff $\set{m\in \omega}{f_{m,j}(n) = \ell} \in D$}.\]
This coincides with the topological $D$-limit. Therefore, we can think of $F:=\set{f_j}{j<i}$ as the $D$-limit of $\seq{F_m}{m<\omega}$, so we define $\lim^D \bar p:=(k,\pi,F)$. Note that $\lim^D \bar p \in Q_{k,\pi,\bar\phi}$.

The sequence $Q_{k,\pi,\bar\phi}$ witness that $\Por^K$ is uniformly $\aleph_0$-uf-$\lim$-linked for any non-principal ultrafilter $D$ on $\omega$. To see this, all that is needed is to prove is~(\ref{Pr_K:eq}) for all $n\in\omega$ where:
\begin{equation}
\parbox{0.8\textwidth}{Given $\bar p^j=\seq{p_m^j}{m\in\omega}\in Q^\omega$ for $j<n$ and $\lim^D\bar p^j\leq q$ for all $j<n$, then $\set{m\in\omega}{q\parallel p_m^j\ \text{for all}\ j<n}\in D$.}
\tag{$\boxplus_n$}
\label{Pr_K:eq}
\end{equation}
To prove~(\ref{Pr_K:eq}), it suffices to show

\begin{clm}
Given $q=(k_q,\pi_q,F_q)\leq \lim^D\bar p$ we have\/ $\set{m\in\omega}{q\parallel p_m}\in D$.
\end{clm}

\begin{varproof}
Since $q=(k_q,\sigma_q,F_q)\leq\lim^D\bar p$, 
either $\sigma_q(f{\upharpoonright}n)=f(n)$ or $\sigma_q(f{\upharpoonright}n+1)=f(n+1)$ for all $f\in F$ and for all $k<n<k_{q-1}$. Then 
\[a_n=\set{m\in\omega}{\text{either $\sigma_q(f_j{\upharpoonright}n)=f(n)$ or $\sigma_q(f_{m}{\upharpoonright}n+1)=f_{m}(n+1)$}}\in D.\]
Hence, $\bigcap\set{a_n}{k<n<k_{q-1}}\neq\emptyset$, so choose $m\in a_n$ for all $k<n<k_{q-1}$ and put $r=(k,\sigma,F_q\cup F_m)$. This is a condition of $\Por^K$. Furthermore, $r$ is stronger than $q$ and $q_m$, so we are done.
\qedsymbol$_\text{Claim}$
\end{varproof}
\let\qed\relax
\end{varproof}

Below, we present Brendle's forcing that increases $\efrak^\const$:

\begin{definition}[{\cite{BreIII}}]\label{b0}
Let $\Por^\omega$ be the poset whose conditions are triples $(k,\sigma, F)$ such that $k\in\omega$, $\sigma\colon{}^{<\omega}\omega\to\omega$ is finite partial function and $F\subseteq\omega^\omega$ is finite, and such that the following requirements are met:
\begin{itemize}
\item
$|s|\leq k$ for all $s\in\dom(\sigma)$,

\item
$f{\upharpoonright}n\in\dom(\sigma)$ for all $f\in F$ and $n\leq k$,
\item
$f{\upharpoonright}k\neq g{\upharpoonright}k$ for all $f\neq g$ belonging to $F$, and
\item
$\sigma(f{\upharpoonright}k)=f(k)$ for all $f\in F$.
\end{itemize}
We order $\Por^\omega$ by $(\ell,\tau, H)\leq(k,\sigma, F)$ iff $\ell\geq k$, $\tau\subseteq\sigma$, $H\supseteq F$, and for all $f\in F$ and for all $k<n<\ell-1$, either $\tau(f{\upharpoonright}n)=f(n)$ or $\tau(f{\upharpoonright}n+1)=f(n+1)$.

Let $G$ be a $\Por^\omega$-generic over $V$. Then in $V[G]$, the generic predictor is given by
\[\sigma_\gen=\bigcup\set{\sigma}{\exists k, F\colon(k,\sigma, F)\in G}.\]
Then it ca be proved that $\sigma_\gen$ is a function from ${}^{<\omega}\omega$ to $\omega$ that predicts constantly $f$, for any $f\in \baire\cap V$. Consequently, $\Por^\omega$ increases $\efrak^\const$. 
\end{definition}

Notice that 

\begin{lemma}[{\cite[Lem.~3.2]{BreIII}}]\label{u1}
$\Por^\omega$ is $\sigma$-linked. Even more, $\Por^\omega$ is $\sigma$-3-linked but it can not be $\sigma$-4-linked (see~\cite{BreShevaIV}).
\end{lemma}

We inspect the proof that $\Por^\omega$ does not add a dominating real to see that $\sigma$-uf-linked (see~\cite[Lem.~3.3]{BreIII}). 

\begin{lemma}\label{t4}
$\Por^\omega$ is $\sigma$-uf-linked.
\end{lemma}

\begin{varproof}
For $k$, $\sigma$ as in~\autoref{b0}, and $\bar\phi=\{\phi_0,\phi_1,\ldots,\phi_{i-1}\}\subseteq\omega^{k}$ let
\[Q_{k,\sigma,\bar\phi}=\set{(k',\sigma',F)\in\Por^\omega}{k=k', \sigma'=\sigma, |F|=i,\setand\forall f\in F\ \exists j<i\colon f{\upharpoonright}k=\phi_j}.\]
We shall prove that $Q_{k,\sigma,\bar\phi}$ is $\Fr$-linked. Let $\bar{p}=\seq{p_m}{m\in\omega}$ be a 
sequence in $Q_{k,\sigma,\bar\phi}$. Write $p_m = (k,\sigma,F_m)$ with $F_m=\set{f_{m,j}}{j<i}$ where $f_{m,j}{\upharpoonright}k=\phi_j$. Using a standard compactness argument, we may assume without loss that for all $j<i$, either
\begin{enumerate}[label=\rm($\boxplus_j$)]
\item\label{b4:I}
there is $g_j\in\baire$ such that $\lim_{m\to\infty}f_{m,j}=g_j$, or
\end{enumerate}
\begin{enumerate}[label=\rm($\boxtimes_j$)]
\item\label{b4:II}
there are $\ell_j\geq k$ and $\psi_j\in\omega^{\ell_j}$ such that $f_{m,j}{\upharpoonright}\ell_j=\psi_j$ for all $m$
and the values $f_{n,m}(\ell_j)$ are all distint. 
\end{enumerate}
Then define $p:= (\ell,\tau,F)\in\Por^\omega$ as follows:
\begin{itemize}
\item
let $A=\set{j<i}{j\text{ fulfills~\ref{b4:II}}}$, so 
for $j\in A$ choose $f^*_j\supseteq \psi_j$ arbitrarily. Put
$F=\set{f^*_j}{j\in A}$;
\item
let $\ell>\ell_j$ for all $j\in A$ and;
\item
extend $\sigma$ to $\tau$ so that $\tau(f^*_j{\upharpoonright}n)=f^*_j(n)$ for all $j$ and all $n$ with $k<n\leq\ell$.
\end{itemize}
We claim that $p\Vdash\lqq\exists^\infty m\in\omega\colon p_m\in\dot G\rqq$. To see this, let $q=(k^q,\sigma^q,F^q)\leq p$ and let $m<\omega$. Next choose $m_0$ such that 
\begin{itemize}
\item
$m_0>m$,
\item
$f_{m,j}{\upharpoonright}k^q=g_j{\upharpoonright}k^q$ for all $j$ satisfying~\ref{b4:I}
\item
$f_{m_0,j}{\upharpoonright}n\notin\dom(\sigma^q)$ for all $j$ satisfying~\ref{b4:II} and $\ell_j<n\leq\ell$, 
\end{itemize}
We then define $r=(k^q,\tau_0,F_{m_0})$ as follows: enlarge $\tau^q$ to $\tau_0$ so that for all $j\in A$ and for all $n$ with $\ell_j< n\leq k^q$, $f_{m,j}{\upharpoonright}n\in\dom(\tau_0)$ and $\tau_0(f_{m,j}{\upharpoonright}n)=f_{m,j}(n)$. It is not hard to see that $r\in\Por^\omega$ and $r\leq (k,\sigma,F_{m_0})$. Even more, $(k^q,\tau_0,F^q)\leq (k^q,\sigma^q,F^q)$. By using~\autoref{u1}, we obtain that $(k^q,\tau_0,F^q)$ and $r$ are compatible, then we are done.
\end{varproof}

\subsection{Simple matrix iterations}
\

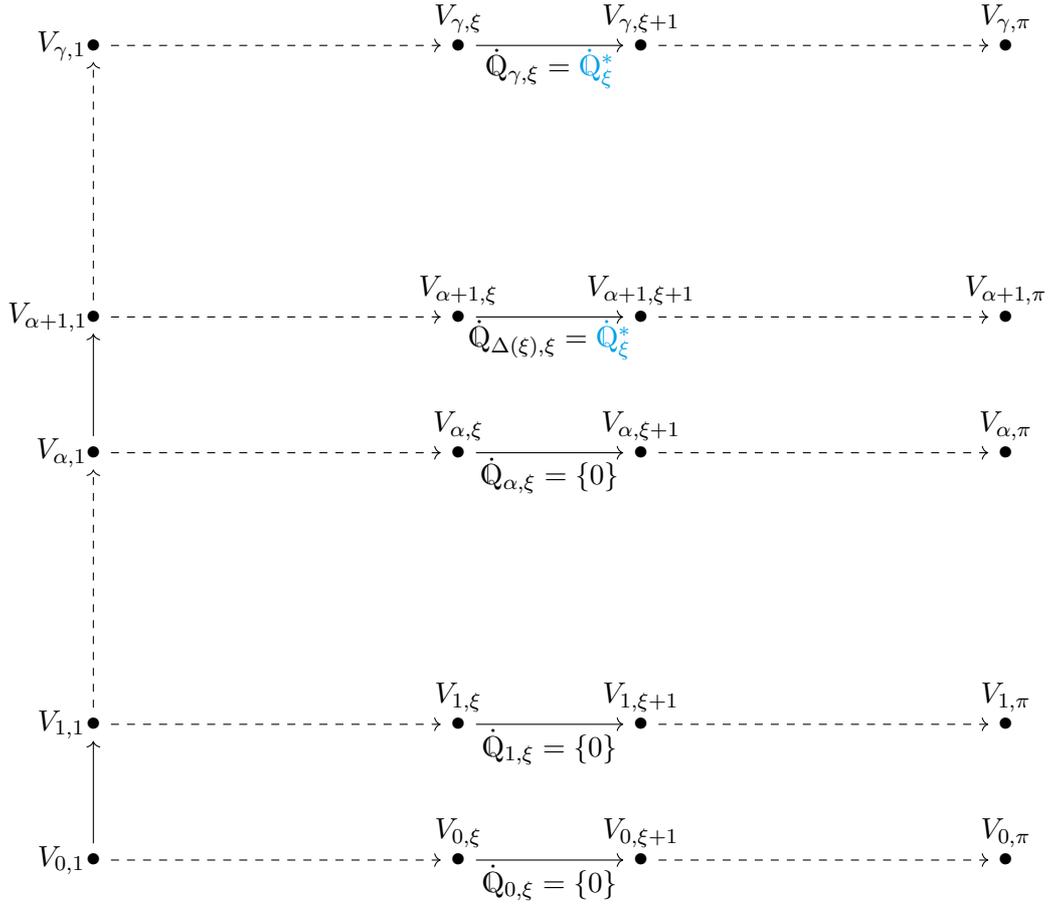
\begin{figure}[ht]
\centering
\begin{tikzpicture}[scale=1.2]
\small{

\node (f00) at (0,0){$\bullet$};
\node (f01) at (0,1.5){$\bullet$};
\node (f02) at (0,4.5) {$\bullet$} ;
\node (f03) at (0,6) {$\bullet$} ;
\node (f04) at (0,9) {$\bullet$} ;
\node (fxi0) at (4,0){$\bullet$};
\node (fxi1) at (4,1.5){$\bullet$};
\node (fxi2) at (4,4.5) {$\bullet$};
\node (fxi3) at (4,6) {$\bullet$} ;
\node (fxi4) at (4,9) {$\bullet$} ;
\node (fxi+10) at (6,0){$\bullet$};
\node (fxi+11) at (6,1.5){$\bullet$};
\node (fxi+12) at (6,4.5) {$\bullet$};
\node (fxi+13) at (6,6) {$\bullet$} ;
\node (fxi+14) at (6,9) {$\bullet$} ;
\node (fpi0) at (10,0){$\bullet$};
\node (fpi1) at (10,1.5){$\bullet$};
\node (fpi2) at (10,4.5) {$\bullet$};
\node (fpi3) at (10,6) {$\bullet$} ;
\node (fpi4) at (10,9) {$\bullet$} ;

\draw   (f00) edge [->] (f01);
\draw[dashed]    (f01) edge [->] (f02);
\draw       (f02) edge [->] (f03);
\draw[dashed]       (f03) edge [->] (f04);

\draw[dashed]   (f00) edge [->] (fxi0);
\draw      (fxi0) edge [->] (fxi+10);
\draw[dashed]       (fxi+10) edge [->] (fpi0);

\draw[dashed]   (f01) edge [->] (fxi1);
\draw        (fxi1) edge [->] (fxi+11);
\draw[dashed]        (fxi+11) edge [->] (fpi1);

\draw[dashed]  (f02) edge [->] (fxi2);
\draw        (fxi2) edge [->] (fxi+12);
\draw[dashed]        (fxi+12) edge [->] (fpi2);

\draw[dashed]  (f03) edge [->] (fxi3);
\draw       (fxi3) edge [->] (fxi+13);
\draw[dashed]      (fxi+13) edge [->] (fpi3);

\draw[dashed]  (f04) edge [->] (fxi4);
\draw        (fxi4) edge [->] (fxi+14);
\draw[dashed]        (fxi+14) edge [->] (fpi4);

\node at (-0.35,0) {$V_{0,1}$};
\node at (-0.35,1.5) {$V_{1,1}$};
\node at (-0.35,4.5) {$V_{\alpha,1}$};
\node at (-0.5,6) {$V_{\alpha+1,1}$};
\node at (-0.35,9) {$V_{\gamma,1}$};

\node at (4,0.3) {$V_{0,\xi}$};
\node at (6,0.3) {$V_{0,\xi+1}$};

\node at (4,1.8) {$V_{1,\xi}$};
\node at (6,1.8) {$V_{1,\xi+1}$};

\node at (4,4.8) {$V_{\alpha,\xi}$};
\node at (6,4.8) {$V_{\alpha,\xi+1}$};

\node at (4,6.3) {$V_{\alpha+1,\xi}$};
\node at (6,6.3) {$V_{\alpha+1,\xi+1}$};

\node at (4,9.3) {$V_{\gamma,\xi}$};
\node at (6,9.3) {$V_{\gamma,\xi+1}$};

\node at (10,0.3) {$V_{0,\pi}$};
\node at (10,1.8) {$V_{1,\pi}$};
\node at (10,4.8) {$V_{\alpha,\pi}$};
\node at (10,6.3) {$V_{\alpha+1,\pi}$};
\node at (10,9.3) {$V_{\gamma,\pi}$};

\node at (5,-0.25) {$\Qnm_{0,\xi}=\{0\}$};
\node at (5,1.25) {$\Qnm_{1,\xi}=\{0\}$};
\node at (5,4.25) {$\Qnm_{\alpha,\xi}=\{0\}$};
\node at (5,5.75) {$\Qnm_{\Delta(\xi),\xi}=\subiii{\Qnm^*_\xi}$};
\node at (5,8.75) {$\Qnm_{\gamma,\xi}=\subiii{\Qnm^*_\xi}$};



}
\end{tikzpicture}
\caption{A simple matrix iteration}
\label{matrixuf}
\end{figure}

\begin{definition}[{\cite[Def.~2.10]{BCM}}]\label{f1}
A~\emph{simple matrix iteration} of ccc posets (\autoref{matrixuf}) is composed of the following objects:
\begin{enumerate}[label=\rm (\Roman*)]
\item
ordinals $\gamma$ (height) and $\pi$ (length);
\item
a function $\Delta\colon\pi\to\gamma$;
\item
a sequence of posets $\seq{\Por_{\alpha,\xi}}{\alpha\leq \gamma,\ \xi\leq \pi}$ where $\Por_{\alpha,0}$ is the trivial poset for any $\alpha\leq \gamma$;
\item
for each $\xi<\pi$, $\Qnm^*_\xi$ is a
$\Por_{\Delta(\xi),\xi}$-name of a~poset such that $\Por_{\gamma,\xi}$ forces it to be ccc;
\item
$\Por_{\alpha,\xi+1}=\Por_{\alpha,\xi}\ast\Qnm_{\alpha,\xi}$, where
\[\Qnm_{\alpha,\xi}:=
\begin{cases}
\Qnm^*_\xi,&\text{if $\alpha\geq\Delta(\xi)$,}\\
\{0\},&\text{otherwise;}
\end{cases}\]
\item
for $\xi$ limit, $\Por_{\alpha,\xi}:=\limdir_{\eta<\xi}\Por_{\alpha,\eta}$.
\end{enumerate}

It is known that $\alpha\leq\beta\leq\gamma$ and $\xi\leq\eta\leq\pi$ imply $\Por_{\alpha,\xi}\subsetdot\Por_{\beta,\eta}$, see e.g.~\cite{B1S} and \cite[Cor.~4.31]{CM}. If $G$ is $\Por_{\gamma,\pi}$-generic
over $V$, we denote $V_{\alpha,\xi}= [G\cap\Por_{\alpha,\xi}]$ for all $\alpha\leq\gamma$ and $\xi\leq\pi$.
\end{definition}

\begin{lemma}[{\cite[Lemma~5]{BrF}, see also~\cite[Cor.~2.6]{mejiamatrix}}]\label{realint}
Assume that\/ $\Por_{\gamma, \pi}$ is a simple matrix iteration as in~\autoref{f1} with\/ $\cf(\gamma)>\omega$.
Then, for any $\xi\leq\pi$,
\begin{enumerate}[label=\rm (\alph*)]
\item
$\Por_{\gamma,\xi}$ is the direct limit of\/ $\seq{\Por_{\alpha,\xi}}{\alpha<\gamma}$, and
\item
if $\eta<\cf(\gamma)$ and $\dot{f}$ is a\/ $\Por_{\gamma,\xi}$-name of a function from $\eta$ into\/ $\bigcup_{\alpha<\gamma}V_{\alpha,\xi}$ then $\dot{f}$ is forced to be equal to a\/ $\Por_{\alpha,\xi}$-name for some $\alpha<\gamma$. In particular, the reals in $V_{\gamma,\xi}$ are precisely the reals in\/ $\bigcup_{\alpha<\gamma}V_{\alpha,\xi}$.
\end{enumerate}
\end{lemma}

Using a Polish relational system that is Tukey-equivalent with $\Cv_\Mwf$ (see \autoref{b16}~\ref{b16:5}) we have the following result.

\begin{theorem}[{\cite[Thm.~5.4]{CM}}]\label{matsizebd}
Let\/ $\Por_{\gamma, \pi}$ be a simple matrix iteration as defined in~\autoref{f1}. Assume that, for any $\alpha<\gamma$, there is some $\xi_\alpha<\pi$ such that\/ $\Por_{\alpha+1,\xi_\alpha}$ adds a~Cohen real $\dot{c}_\alpha\in X$ over~$V_{\alpha,\xi_\alpha}$.
Then, for any $\alpha<\gamma$, $\Por_{\alpha+1,\pi}$ forces that $\dot{c}_{\alpha}$ is Cohen over~$V_{\alpha,\pi}$.

In addition, if\/ $\cf(\gamma)>\omega_1$ and $f\colon\cf(\gamma)\to\gamma$ is increasing and cofinal, then\/ $\Por_{\gamma,\pi}$ forces that\/ $\set{\dot{c}_{f(\zeta)}}{\zeta<\cf(\gamma)}$ is a strongly\/ $\cf(\gamma)$-$\Cv_\Mwf$-unbounded family. In particular, $\Por_{\gamma,\pi}$ forces $\gamma\leqT\Cv_\Mwf$ and\/ $\non(\Mwf)\leq\cf(\gamma)\leq\cov(\Mwf)$.
\end{theorem}

The following notion formalizes the matrix iterations with ultrafilters from~\cite{BCM}. 

\begin{definition}[{\cite[Def.~4.2]{BCM}}]\label{f4}
Let $\theta\geq\aleph_1$ and let $\Por_{\gamma, \pi}$ be a simple matrix iteration as in~\autoref{f1}. Say that $\Por_{\gamma, \pi}$ is a~\emph{${<}\theta$-uf-extendable matrix iteration} if for each $\xi<\pi$, $\Por_{\Delta(\xi),\xi}$ forces that $\Qnm_\xi$ is a ${<}\theta$-uf-linked poset.
\end{definition}

\begin{theorem}[{\cite[Thm.~4.4]{BCM}}]\label{mainpres}
Assume that $\theta\leq\mu$ are uncountable cardinals with $\theta$ regular. Let\/ $\Por_{\gamma,\pi}$ be a ${<}\theta$-uf-extendable matrix iteration as in~\autoref{f4} such that
\begin{enumerate}[label=\rm (\roman*)]
\item $\gamma\geq\mu$ and $\pi\geq\mu$,
\item for each $\alpha<\mu$, $\Delta(\alpha)=\alpha+1$ and\/ $\Qnm_\alpha$ is Cohen forcing, and
\item $\dot c_\alpha$ is a\/ $\Por_{\alpha+1,\alpha+1}$-name of the Cohen real in $\omega^\omega$ added by\/ $\Qnm_\alpha$.
\end{enumerate}
Then\/ $\Por_{\alpha,\pi}$ forces that\/ $\set{\dot c_\alpha}{\alpha<\mu}$ is $\theta$-$\baire$-unbounded, in particular, $\Cv_{[\mu]^{<\theta}}\leqT \baire$.
\end{theorem}

\section{Tools for controlling the additivity of~\texorpdfstring{$\Nwf$}{}}\label{sec:s3}

We introduce a new property (which we call $\Pr^2_{\bar n^*}(\lambda)$ property, see~\autoref{def:ppies}) that allows us to maintain the additivity of the null ideal small under certain forcing extensions. Certainly, this resembles a type of property such as the Knaster property. Concretely, we prove~\autoref{Thm:a3} and show in~\autoref{cor:rhoPr2} that $\Pr^2_{\bar n^*}(\lambda)$ is weaker than the $\sigma$-$\bar\rho$-linked property (see~\autoref{d1}). Additionally, we show that various forcing notions fulfill the property $\Pr^2_{\bar n^*}(\lambda)$. 

We start with the definition of the $\sigma$-$\bar\rho$-linked introduced in~\cite{CRS}, which lies $\sigma$-centered and $\sigma$-linked.

\begin{definition}[{\cite[Def.~4.1]{CRS}}]
\label{d1}
Let $\Por$ be a~forcing notion and let $\bar\rho=\seq{\rho_n}{n\in\omega}$ be a sequence of functions $\rho_n\colon\omega\to\omega$ such that $\lim_{k\to\infty}\rho_n(k)=\infty$, $\rho_n(k)\leq k$, $\rho_n(k+1)\geq2$ and $\rho_n\geq\rho_{n+1}$.
\begin{enumerate}[label=\rm(\arabic*)]
\item\label{d1:1}
Say that $\Por$ is \emph{$\sigma$-$\bar\rho$-linked} if there is a sequence $\set{P_n}{n\in\omega}$ such that:
\smallskip
\begin{enumerate}[label=\rm(\alph*)]
\item\label{d1:a}
$\bigcup_{n\in\omega}P_n$ is dense in $\Por$ and
\item\label{d1:b}
for all $n$, $k<\omega$, if $\set{p_i}{i<k}\subseteq P_n$, then there is a subset $Q$ of $\set{p_i}{i<k}$ of size $\rho_n(k)$ that has a common lower bound.
\end{enumerate}

\item
Say that $\Por$ is \emph{$\sigma$-$\bar\rho{\star}$-linked},
if there is a sequence $\seq{P_n}{n\in\omega}$ such that:
\smallskip
\begin{enumerate}[label=\rm(\alph*)]
\item $\Por=\bigcup_{n\in\omega}P_n$ and
\item $\forall n,k\in\omega\
\forall\seq{p_i}{i<k}\in{}^kP_n\
\exists u\in[k]^{\rho_n(k)}$ the set
$\set{p_i}{i\in u}$ has a lower bound.
\end{enumerate}
\item We say that $\Por$ is \emph{$\sigma{\uparrow}$-$\bar\rho$-linked} or \emph{$\sigma{\uparrow}$-$\bar\rho{\star}$-linked}, respectively,
if moreover, $P_n\subseteq P_{n+1}$ for all $n\in\omega$.
\end{enumerate}
\end{definition}

\begin{example}[{\cite[Lem.~4.3]{CRS}}]\label{rho:exam}
\ 
\begin{enumerate}[label=\rm(\arabic*)]
\item\label{exa:1}
Any $\sigma$-centered posets are $\sigma$-$\bar\rho$-linked for $\rho_n(k)=k$.
\item\label{exa:2}
Random forcing\/~$\Bor$ is $\sigma$-$\bar\rho$-linked for some~$\bar\rho$.
\end{enumerate}
\end{example}

Let $\bar\rho=\seq{\rho_n}{n\in\omega}$ as in~\autoref{d1}.
For $h\in\baire$ let
\begin{enumerate}
\item
$g_m(k)=\max\set{n\in\omega}{\rho_m(n)\leq h(k)}$ and $h'(k)=\sum_{m\le k}g_m(k)$,
\item
$\Hwf_{\bar\rho,h}=\set{h_n}{n\in\omega}\subseteq\baire$ where $h_0=h$ and $h_{n+1}=(h_n)'$.
\end{enumerate}
Note that $g_m(k)\geq h(k)$ because $2\leq\rho_m(n)\le n$
(for $n\geq 2$) and hence $h'(k)\geq(k+1)\cdot h(k)$.
Therefore $h_{n+1}(k)\geq(k+1)^{n+1}$ whenever $h(k)\geq1$.

\begin{fact}\label{fct:rho}
Let $\bar\rho=\seq{\rho_n}{n\in\omega}$ as in~\autoref{d1}. Then\/ $\Lc^*_{\Hwf_{\bar\rho,h}}\eqT\Nwf$ whenever
$h$~is positive. Even more, if $h\in\Hwf_{\bar\rho,h}$, then there is an $h^*\in\Hwf_{\bar\rho,h}$ such that $\frac{h}{h^*}$ coverge to $0$.
\end{fact}

\autoref{fct:rho} and ~\autoref{b16}~\ref{b16:3} together yield:

\begin{remark}\label{rem:rancent}
Both subalgebra of random forcing and $\sigma$-centered posets are $\Lc^*_{\Hwf_{\bar\rho,h}}$-good.
\end{remark}

The following result displays that the $\sigma$-$\bar\rho$-linked property behaves well to control $\add(\Nwf)$ small in extension generics. 

\begin{lemma}[{\cite[Cor.~4.16]{CRS}}]\label{goodrho}
Let $\bar\rho$ be as in~\autoref{d1}. Assume
that $h\in\baire$ converges to infinity. Then there is some $\Hwf_{\bar\rho,h}=\set{h_n}{n\in\omega}\subseteq\baire$ with $h_0=h$ such that any $\sigma$-$\bar\rho$-linked posets is\/ $\Lc^*_{\Hwf_{\bar\rho,h}}$-good.
\end{lemma}

Unless otherwise specified, until the end of this section, we fix a poset $\Por$ and a regular uncountable cardinal $\lambda$.

\begin{definition}\label{def:lamrho}
Let $\bar\rho=\seq{\rho_n}{n\in\omega}$ be as in~\autoref{d1}.
\begin{enumerate}
\item
$\Por$ has \emph{$\lambda$-$\bar\rho$-cc} if for every sequence
$\bar p\in{}^\lambda\Por$ for some $A\in[\lambda]^\lambda$ and
for every $m$ large enough, if $B\in[A]^{<\omega}$ then there is
$F\in[B]^{\rho_m(|B|)}$ such that $\set{p_\xi}{\xi\in F}$ has
an lower bound.

\item
$\Por$ has \emph{weak $\lambda$-$\bar\rho$-cc} if for every sequence
$\bar p\in{}^\lambda\Por$ for some $A\in[\lambda]^\lambda$
and $m<\omega$, if $B\in[A]^{<\omega}$ then there is
$F\in[B]^{\rho_m(|B|)}$ such that $\set{p_\xi}{\xi\in F}$ has
an lower bound.

\item
$\Por$ has \emph{strong $\lambda$-$\bar\rho$-cc} if for every sequence
$p\in{}^\lambda\Por$ there is $A\in[\lambda]^\lambda$ such that
for every $B\in[A]^{<\omega}$ and every $m\in\omega$ there is
$F\in[B]^{\rho_m(|B|)}$ such that $\set{p_\xi}{\xi\in F}$ has an
lower bound.
\end{enumerate}
\end{definition}

\begin{definition}
Let $\rho\in\baire$. Assume $\rho(k)\leq k$ and $\lim_k\rho(k)=\infty$.
Say that $\Por$ is $\lambda$-$\rho$-Knaster, if
$\forall p\in{}^\lambda\Por$
$\forall k\in\omega$
$\exists A\in[\lambda]^\lambda$
$\forall u\in[A]^k$
$\exists v\in[u]^{\rho(k)}$
$\set{p_\alpha}{\alpha\in v}$ has an lower bound.
\end{definition}

\begin{fact}\label{fac:z0}
Let $\bar\rho$ be as in~\autoref{d1}. Then 
\begin{enumerate}[label=\rm(\arabic*)]
\item\label{fac:z0:1} Any $\lambda$-$\bar\rho$-cc posets satisfies the weak $\lambda$-$\bar\rho'$-cc for some $\bar\rho'$. 
\item\label{fac:z0:3}  Any strong $\lambda$-$\bar\rho$-cc posets satisfies the $\lambda$-$\rho'$-Knaster property for some $\rho'$.
\end{enumerate}    
\end{fact}

\begin{varproof}
The proof~\ref{fac:z0:1} is clear, so we are left with checking~\ref{fac:z0:1} and~\ref{fac:z0:3}. For~\ref{fac:z0:3}, define $\rho'=\sup\bar\rho$ where $\sup\bar\rho(k)=\max\set{\rho_m(k)}{k\in\omega}$. It is not hard to verify that if $\Por$ has strong $\lambda$-$\bar\rho$-cc, then satisfies the $\lambda$-$\rho'$-Knaster property. Lastly, for~\ref{fac:z0:1} define $\rho=\sup_m\rho_m$, i.e., $\rho(k)=\max\set{\rho_m(k)}{k\in\omega}\leq k$. Note that if $\Por$ has strong $\lambda$-$\bar\rho$-cc, then $\lambda$-$\bar\rho'$-cc for $\rho'=\sup\bar\rho$. 
\end{varproof}

For sequences of functions $\bar\rho$ and~$\bar\sigma$ define
$\bar\rho\bar\sigma=\bar\rho\cdot\bar\sigma=
\seq{\rho_n\circ\sigma_n}{n\in\omega}$.
Denote $d_n(k)=\lfloor k/(n+1)\rfloor$ and $r_n(k)=\lfloor\sqrt k\rfloor$.
Then,
$(\bar\rho\bar d)_n=\rho_n(\lfloor k/(n+1)\rfloor)$,
$(\bar\rho\bar r)_n=\rho_n(\lfloor\sqrt k\rfloor)$, and
$(\bar\rho\bar r\bar d)_n=\rho_n(\lfloor\sqrt{k/(n+1)}\rfloor)$.

\begin{fact}\label{fct}
Assume that\/ $\cf(\lambda)>\omega$.
\begin{enumerate}[label=\rm(\arabic*)]
\item\label{fct-a}
$\sigma$-$\bar\rho{\star}$-linked${}\Rightarrow{}\sigma$-$\bar\rho$-linked${}\Rightarrow{}\sigma$-$\bar\rho\bar r{\star}$-linked.

\item\label{fct-c}
$\sigma$-$\bar\rho{\star}$-linked${}\Rightarrow
\sigma{\uparrow}$-$\bar\rho\bar d{\star}$-linked.

\item\label{fct-d}
$\sigma{\uparrow}$-$\bar\rho{\star}$-linked$
{}\Rightarrow\lambda$-$\bar\rho$-cc. 

\item\label{fct-e}
$\sigma$-$\bar\rho{\star}$-linked${}\Rightarrow\lambda$-$\bar\rho\bar d$-cc.

\item\label{fct-f}
$\sigma$-$\bar\rho$-linked${}\Rightarrow\lambda$-$\bar\rho\bar r\bar d$-cc.
\end{enumerate}
\end{fact}

\begin{varproof}
In~\autoref{d1} we can assume that $\Por=\bigcup_{n\in\omega}P_n$.

\ref{fct-a}
Easy.
To justify the second implication it is enough to realize that every sequence
of length~$k$ either contains $\lfloor\sqrt k\rfloor$ mutually different terms
or $\lfloor\sqrt k\rfloor$ occurrences of the same term.

\ref{fct-c}
Define $P'_n=\bigcup_{m\leq n}P_m$ for $n\in\omega$.
Every sequence of length~$k$ in $P'_n$ has a~subsequence of length
$\lfloor k/(n+1)\rfloor$ in some $P_m$ with $m\leq n$.

\ref{fct-e}~is a~consequence of \ref{fct-c} and \ref{fct-d}.

\ref{fct-f}~is a~consequence of \ref{fct-a} and \ref{fct-e}.
\end{varproof}

\begin{definition}\label{def:ppies}
Let $\bar n^*=\seq{n_i}{i\in\omega}\in\baire$.
\begin{enumerate}[label=\rm(\arabic*)]
\item\label{def:ppies:a}
Say that $\Por$ satisfies the $\Pr^1(\lambda)$ property if
for every $\bar p\in{}^{\lambda}\Por$ there are some
$S\in[\lambda]^\lambda$ and $k\geq 2$ such that
$\forall u\in[S]^{\leq k}$ $\set{p_\alpha}{\alpha\in u}$ has
a~common lower bound.
\item[--]
$\Por$~satisfies $\Pr^1(\lambda)\Leftrightarrow
\forall\bar p\in{}^\lambda\Por\
\exists S\in[\lambda]^\lambda\
\forall\alpha,\beta\in S$
$p_\alpha$ and $p_\beta$ are compatible.

\item\label{def:ppies:b}
Say that $\Por$ satisfies the property $\Pr^2_{\bar n}(\lambda)$
if for every $\bar p\in{}^{\lambda}\Por$ there are some
$S\in[\lambda]^\lambda$ and $k\geq 2$ such that
$\forall i>k$
$\forall u\in[S]^{n_{i+k}}$
$\exists v\in[u]^{n_i}$
$\set{p_\alpha}{\alpha\in v}$
has a~common lower bound.
\item Say that $\Por$ satisfies the strong property $\Pr^2_{\bar n}(\lambda)$ if
$\forall\bar p\in{}^\lambda\Por$
$\exists S\in[\lambda]^\lambda$
$\forall^\infty i\in\omega$
$\forall u\in[S]^{n_{i+1}}$
$\exists v\in[u]^{n_i}$
$\set{p_\alpha}{\alpha\in v}$
has a~common lower bound.
\end{enumerate}
\end{definition}

\begin{remark}
Observe that Knaster property is exactly the $\Pr^1(\omega_1)$ property.
\end{remark}

\begin{lemma}\label{lem:ppies}
Let $\bar n\in\baireinc$.
\begin{enumerate}[label=\rm(\arabic*)]
\item\label{lem:ppies:a}
$\Pr^2_{\bar n}(\lambda)\Rightarrow{}$strong property $\Pr^2_{\bar m}(\lambda)$ where $m_i=n_{2^i}$.
\item\label{lem:ppies:b}
Strong property $\Pr^2_{\bar n}(\lambda)\Rightarrow\Pr^2_{\bar n}(\lambda)$.
\end{enumerate}
\end{lemma}

\begin{varproof}
Let $\bar p\in{}^\lambda\Por$ and $S\in[\lambda]^\lambda$.

\ref{lem:ppies:a}: Assume that $k\ge2$ is such that
$\forall u\in[S]^{\leq k}$ $\set{p_\alpha}{\alpha\in u}$ has
a~common lower bound.
Then for every $i$ with $2^i>k$ we have $m_{i+1}\geq n_{2^i+k}$
and hence for every $u\in[S]^{m_{i+1}}$ there is
$v\in[u]^{n_{2^i}}=[u]^{m_i}$ such that $\set{p_\alpha}{\alpha\in v}$
has a~common lower bound.

\ref{lem:ppies:b}: Assume that $k\geq 2$ and for every $i>k$ and $u\in[S]^{n_{i+1}}$,
$\exists v\in[u]^{n_i}$ $\set{p_\alpha}{\alpha\in v}$
has a~common lower bound.
Then the same holds also for every $u\in[S]^{n_{i+k}}$
and hence it holds that
$\forall u\in[S]^{\leq k}$ $\set{p_\alpha}{\alpha\in u}$ has
a~common lower bound.
\end{varproof}

\begin{lemma}\label{lem:ppties}
Let $\bar\rho$ be as in~\autoref{d1}.
\begin{enumerate}[label=\rm(\arabic*)]
\item\label{lem:ppties:a}
If\/ $\Por$ is $\sigma$-$\bar\rho$-linked, then\/ $\Por$~satisfies\/
strong property $\Pr^2_{\bar n}(\lambda)$ for some $\bar n\in\incrbaire$.
\item\label{lem:ppties:b}
Any weak $\lambda$-$\bar\rho$-cc posets satisfies the strong property $\Pr^2_{\bar n}(\lambda)$ for some $\bar n\in\incrbaire$
\end{enumerate}
\end{lemma}

\begin{varproof}
\ref{lem:ppies:a}: Let $\Por=\bigcup_{n\in\omega}P_n$ be such that
for all $n$, $k<\omega$, if $\set{p_i}{i<k}\subseteq P_n$,
then there is a~subset of $\set{p_i}{i<k}$ of size $\rho_n(k)$
that has a~common lower bound.
Let $f\in\baire$ be a~strictly increasing function such that
$f(m)\geq\min\set{k\in\omega}{\forall n\leq m\bsp\forall k'\geq k\bsp\rho_n(k')\geq m}$.
By induction define $n_0=0$ and $n_{i+1}=f(n_i)$.
Assume $\bar p\in{}^\lambda\Por$.
Since $\lambda$~is regular there is $n\in\omega$ such that
the set $S=\set{\alpha<\lambda}{p_\alpha\in P_n}$ has
cardinality~$\lambda$.
Then for all but finitely many~$i$, $n_i\geq n$ and for every such~$i$,
$\rho_n(n_{i+1})=\rho_n(f(n_i))\geq n_i$.
This implies the property strong property $\Pr^2_{\bar n}(\lambda)$.

\ref{lem:ppies:b}: Define $f$ and $\bar n$ as in~\ref{lem:ppies:a}.
Assume $\bar p\in{}^\lambda\Por$ and let $A\in[\lambda]^\lambda$ and $n<\omega$ be such that
$\forall B\in[A]^{<\omega}$ $\exists F\in[B]^{\rho_n(|B|)}$ such that
$\set{p_\xi}{\xi\in F}$ has an lower bound.
If $n_i\geq n$, then $\rho_n(n_{i+1})=\rho_n(f(n_i))\geq n_i$ and hence
for every $u\in[A]^{n_{i+1}}$ there is $v\in[u]^{n_i}$ such that
$\set{p_\alpha}{\alpha\in v}$
has a~common lower bound.
It follows that the strong property property $\Pr^2_{\bar n}(\lambda)$ holds.
\end{varproof}

As a direct consequence of~\autoref{lem:ppies}~\ref{lem:ppies:b}, we obtain:

\begin{corollary}\label{cor:rhoPr2}
Any $\sigma$-$\bar\rho$-linked satisfies the property\/ $\Pr^2_{\bar n}(\lambda)$. 
\end{corollary}

Using~\autoref{lem:ppties} and~\autoref{rho:exam}, we can infer:

\begin{example}\label{exmp:Pr^2}
\
\begin{enumerate}
\item
Any $\sigma$-centered posets satisfies $\Pr^2_{\bar n}(\lambda)$ property.
\item
Random forcing satisfies $\Pr^2_{\bar n}(\lambda)$ property.
\end{enumerate}
\end{example}

We now aim to prove that $\Pr^2_{\bar n^*}(\lambda)$ is preserved in FS iterations of ccc posets.

\begin{lemma}\label{lem:stsuc}
Assume that $\lambda$~is a~regular cardinal number.
If\/ $\Por$~satisfies\/ $\Pr^2_{\bar n}(\lambda)$ and\/
$\Vdash_\Por\lqq\Qnm$~satisfies\/ $\Pr^2_{\bar n}(\lambda)\rqq$,
then\/
$\Por*\Qnm$ satisfies\/ $\Pr^2_{\bar n}(\lambda)$.  
\end{lemma}

\begin{varproof}
Let $\seq{(p_\xi,\dot q_\xi)}{\xi<\lambda}\in{}^{\lambda}(\Por*\dot \Qor)$.
Since $\Por$ has $\lambda$-cc and $\lambda$~is regular, there is $\xi_0<\lambda$ such that
for every $\xi\geq\xi_0$
the sets $\set{p_\zeta}{\zeta\geq\xi}$ are pairwise \lqq forcing equivalent\rqq, i.e.,
$\forall\xi,\zeta\geq\xi_0\bsp\forall r\leq p_\xi\bsp\exists r'\leq r\bsp
\exists\zeta'\geq\zeta\bsp r'\leq p_{\zeta'}$.
Let $\dot A_1$ denote the $\Por$-name
$\set{\zeta\geq\xi_0}{p_\zeta\in\dot G}$
where $\dot G$ is the canonical $\Por$-name for a~generic subset
of~$\Por$.
Then for every $\xi\geq\xi_0$, $p_\xi\Vdash\lqq\dot A_1\in[\lambda]^\lambda\,\rqq$.
As $\Vdash_\Por\lqq\Qnm$~satisfies\/ $\Pr^2_{\bar n}(\lambda)\rqq$, 
there are $\Por$-names
$\dot A_2$ and $\dot k\geq2$ such that for every $\xi\geq\xi_0$,
$p_\xi$~forces
\begin{enumerate}
\item[$(\boxdot_1)$]
$\dot A_2\in[\dot A_1]^\lambda$ and
$\forall i>\dot k$
$\forall u\in[\dot A_2]^{n_{i+\dot k}}$
$\exists v\in[v]^{n_i}$
$\set{\dot q_\zeta}{\zeta\in v}$
has a~common lower bound.
\end{enumerate}
Choose $p_*\leq p_{\xi_0}$ such that for some $k_0\in\omega$,
$p_*\Vdash\dot k=k_0$.
Let
\[
A_3=\set{\zeta\geq\xi_0}
{\exists p_\zeta'\leq p_*\bsp p_\zeta'\leq p_\zeta\setand
p_\zeta'\Vdash\lqq\zeta\in\dot A_2\rqq}.
\]
Then $|A_3|=\lambda$ because
$p_*\Vdash\lqq\dot A_1\in[\lambda]^\lambda\setand(\boxdot_1)\setand\dot A_2\subseteq A_3\rqq$.
For each $\zeta\in A_3$ choose $p_\zeta'\in\Por$ as in the
definition of~$A_3$.
Applying the fact that $\Por$~satisfies\/ $\Pr^2_{\bar n^*}(\lambda)$ to the set
$\seq{p'_\zeta}{\zeta\in A_3}$ we find $S\in[A_3]^\lambda$ and $k_1\geq2$ such that
\begin{enumerate}
\item[$(\boxdot_2)$]
$\forall i\geq k_1$
$\forall u\in[S]^{n_{i+k_1}}$
$\exists v\in[u]^{n_i}$
the set $\set{p_\zeta'}{\zeta\in v}$ has a~common lower bound.
\end{enumerate}

We verify that $S$ and $k=k_0+k_1$ witness $\Pr^2_{\bar n}(\lambda)$
for $\Por*\Qnm$.
Let $i>k_0+k_1$ and $u\in[A_4]^{n_{i+k_0+k_1}}$.
Since $i+k_0>k_1$, by~$(\boxdot_2)$ there is
$v_0\in[u]^{n_{i+k_0}}$ such that
$\set{p_\zeta'}{\zeta\in v_0}$ has a~common lower bound~$r$
and by definition of~$A_3$, $r\Vdash\lqq v_0\in[\dot A_2]^{n_{i+k_0}}\rqq$.
Since $i>k_0$ and $r\Vdash\lqq\dot k=k_0\setand(\boxdot_1)\rqq$ there is $v_1\in[v_0]^{n_i}$
such that some extension of~$r$
(which is a~common extension of $\set{p_\zeta}{\zeta\in v_0}$)
forces that $\lqq\set{\dot q_\zeta}{\zeta\in v_1}$ has
a~common lower bound\rqq.
It follows that $\set{(p_\zeta,\dot q_\zeta)}{\zeta\in v_1}$
has a~common lower bound in $\Por*\Qnm$.
\end{varproof}

\begin{corollary}\label{Pr^2:preservation}
Any FS iteration of posets satisfying\/ $\Pr^2_{\bar n^*}(\lambda)$
satisfies\/ $\Pr^2_{\bar n}(\lambda)$.
\end{corollary}
\begin{varproof}
Let $\seq{\Por_\alpha,\Qnm_\beta}{\beta<\pi\setand\alpha\leq\pi}$
be a~FS iteration of posets such that for every $\alpha<\pi$,
$\Vdash_\beta\lqq\Qnm_\beta$~satisfies~$\Pr^2_{\bar n}(\lambda)\rqq$.
By induction on $\alpha\leq\pi$ we prove that
$\Por_\alpha$~satisfies~$\Pr^2_{\bar n}(\lambda)$.
The case $\alpha=0$ is trivial and the successor step folows by \autoref{lem:stsuc}.
Assume that $\alpha$~is a~limit ordinal and let $\seq{p_\xi}{\xi<\lambda}\in{}^\lambda\Por_\alpha$.
By the $\Delta$-system lemma there is $\beta<\alpha$ and $S_0\in[\lambda]^\lambda$ such that
\[
\tag{$\star$}\label{Delta-system}
\forall\xi,\zeta\in S_0\bsp
(\xi\neq\zeta\Rightarrow
(\dom(p_\xi)\cap\dom(p_\zeta))\smallsetminus\beta=\emptyset).
\]
By the induction hypothesis $P_\beta$~satisfies $\Pr^2_{\bar n}(\lambda)$ and hence
there is $S\in[S_0]^\lambda$ and $k\geq2$ such that
for every $i\geq k$ and $u\in[S]^{n_{i+k}}$ there is $v\in[u]^{n_i}$ such that
\[
\text{$\set{p_\zeta{\restriction}\beta}{\zeta\in v}$ has a~common lower bound in~$\Por_\beta$}
\]
and then, due to~\eqref{Delta-system}, also $\set{p_\zeta}{\zeta\in v}$ has a~common lower bound in~$\Por_\alpha$.
\end{varproof}

\autoref{def:ppies}~\ref{def:ppies:a} is weaker than $\sigma$-linked:

\begin{lemma}\label{Pr1-linked}
If\/ $\Por$ is $\sigma$-linked and $\cf(\lambda)>\omega$, then\/ $\Por$~has the\/ $\Pr^1(\lambda)$ property.
\qed
\end{lemma}

\begin{lemma}\label{Pr-ppty}
Let\/ $2\leq K<\omega$. Then 
\begin{enumerate}[label=\rm(\arabic*)]
\item\label{Pr-ppty:a}
$\Por^K$ satisfies the\/ $\Pr^2_{\bar n}(\lambda)$ property when $n_{i+1}\geq K^{n_i}$.
\item\label{Pr-ppty:b}
$\Por^K$ is $\sigma$-$\bar\rho$-linked for some $\bar\rho$.
\end{enumerate}
\end{lemma}

\begin{varproof}
For every $p\in\Por^K$ consider the enumeration
$\seq{\eta_p^i}{i<|F_p|}$ of~$F_p$ which agrees with the
lexicographic arrangement and let
$\bar s_p=
\seq{\eta_p^i{\restriction}n_p}{i<|F_p|}\in
{}^{|F_p|}(^{n_p}K)$.
Denote
\begin{align*}
Q&=\set{(n_p,\pi_p,|F_p|,\bar s_p)}{p\in\Por^K},\\
P_{n,\pi,k,\bar s}&=
\set{p\in\Por^K}{(n_p,\pi_p,|F_p|,\bar s_p)=(n,\pi,k,\bar s)},\quad
(n,\pi,m,\bar s)\in Q.
\end{align*}
Clearly, $Q$~is a~countable set and
$\Por^K=\bigcup_{(n,\pi,k,\bar s)\in Q}P_{n,\pi,k,\bar s}$.
Every set $P_{n,\pi,k,\bar s}$ is $2$-linked
(see the proof of this fact in~\cite[Lem.~3.2]{BreIII} in the case
of~$\Por^\omega$).

By induction on~$k$ define $e_k(m)$ by
$e_0(m)=m$ and $e_{k+1}(m)=K^{e_k(m)}$.
Let $\bar p=\set{p_\alpha}{\alpha<\lambda}\in{}^\lambda(\Por^K)$
in case~\ref{Pr-ppty:a}.
Since $Q$~is countable there is $(n,\pi,k,\bar s)\in Q$ such that
the set $S=\set{\alpha<\lambda}{p_\alpha\in P_{n,\pi,k,\bar s}}$
has cardinality~$\lambda$.
In case~\ref{Pr-ppty:b} define
$\rho_{n,\pi,k,\bar s}(m)=2$, if $2\leq m<n_k$, and
$\rho_{n,\pi,k,\bar s}(m)=n_i$ if $n_{i+k}\leq m<n_{i+k+1}$.
As $n_{i+k}\geq e_k(n_i)$, in both cases it is enough to prove
\begin{enumerate}[label=\rm($\boxdot$)]
\item\label{lem:eqp}
$\forall \bar p\in{}^{e_k(n_i)}P_{n,\pi,k,\bar s}\ 
\exists v\in[e_k(n_i)]^{n_i}\ 
\set{p_j}{j\in v}$ has a lower bound.
\end{enumerate}

For different $\eta,\eta'\in{}^\omega K$ let
$\nu(\eta,\eta')=\min\set{i\in\omega}{\eta(i)\neq\eta'(i)}$.
If $p,p'\in P_{n,\pi,k,\bar s}$
(where $\bar s=\seq{s_i}{i<k}\in{}^k(^n K)$)
and $i<k$,
then $\eta_p^i,\eta_{p'}^i\in[s_i]\cap{}^\omega K$ and therefore
$\eta_p^i\neq\eta_{p'}^i\Rightarrow\nu(\eta_p^i,\eta_{p'}^i)\geq n$.

\begin{clm}\label{K-clm}
Assume that $\bar p=\seq{p_j}{j<\ell}\in{}^\ell P_{n,\pi,k,\bar s}$
is a~sequence of conditions such that for every $i<k$ and every
$j_0<j_1<j_2<\ell$,
\[
\text{if $\eta_{p_{j_0}}^i\neq\eta_{p_{j_1}}^i$ and
$\eta_{p_{j_1}}^i\neq\eta_{p_{j_2}}^i$, then
$\eta_{p_{j_0}}^i\neq\eta_{p_{j_2}}^i$ and
$\nu(\eta_{p_{j_0}}^i,\eta_{p_{j_1}}^i)<\nu(\eta_{p_{j_1}}^i,\eta_{p_{j_2}}^i)$}.
\]
Then there is a~common extension of the conditions in $\bar p$.
\end{clm}

\begin{varproof}
Let $n'\geq n$ be the least natural number such that for every $i<k$,
$\eta_{p_j}^i{\restriction}n'$ for $j<\ell$ are pairwise different, i.e.,
$n'=\max(\{n\}\cup\set
{1+\nu(\eta_{p_j}^i,\eta_{p_{\ell-1}}^i)}
{i<k\comma j<\ell\comma\eta_{p_j}\neq\eta_{\ell-1}})$.
Let $F'=\bigcup_{j<\ell}F_{p_j}$
and define $\pi'\to{}^{\leq n'}K\to K$ so that $\pi\subseteq\pi'$ and
for every $s\in{}^{\leq n'}K\smallsetminus{}^{\leq n}K$,
\[
\pi'(s)=
\begin{cases}
\eta_{p_j}^i(|s|),&
\text{if $j<\ell$ is maximal such that $s_i\subseteq s\subseteq\eta_{p_{j}}^i$},\\
0,&\text{if there are no $i<k$ and $j<\ell$ such that $s_i\subseteq s\subseteq\eta_{p_{j}}^i$.}
\end{cases}
\]
Now, $(n',\pi',F')\in\Por^K$ and for every $j<\ell$,
$(n',\pi',F')$ is an extension of~$p_j$
because for every $i<k$ there is at most one~$x\in(n,n')$
such that
$\pi'(\eta_{p_j}^i{\restriction}x)\neq
\eta_{p_j}^i(x)$.
\end{varproof}

We continue with the proof of~\autoref{Pr-ppty} by
the following simple observation:
\begin{enumerate}[label=\rm($\boxast$)]
\item\label{eq:Pr}
If $a$ is a~finite set, $s\in{}^{<\omega}K$, and
$\bar\eta=\seq{\eta_p}{p\in a}\in{}^a([s]\cap{}^\omega K)$, then
there is $p_0\in a$ such that either
$\exists s'\supsetneq s$
$\exists a'\in[a]^{\geq|a|/K}$
$\forall p\in a'$
$\eta_p\in[s']\setand\eta_{p_0}\notin[s']$, or
$\forall p\in a$ $\eta_p=\eta_{p_0}$.
\end{enumerate}
To see this, assume that $\eta_p\neq\eta_{p'}$ for some $p,p'\in a$.
Denote
$x=\min
\set{\nu(\eta_p,\eta_{p'})}{\eta_p\neq\eta_{p'}\setand p,p'\in a}$
and find $y_0\neq y_1$ in~$K$ such that
$\set{p\in a}{\eta_p(x)=y_0}\neq\emptyset$ and
$|\set{p\in a}{\eta_p(x)=y_1}|\geq|a|/K$.
Now choose $p_0\in a$ such that $\eta_{p_0}(x)=y_0$ and let
$a'=\set{p\in a}{\eta_p(x)=y_1}$,

By repeating~\ref{eq:Pr} $\ell$~times, if $|a|\geq K^\ell$, we
obtain a~one-to-one sequence $\seq{p_j}{j<\ell}\in{}^{\ell}a$
such that for every $j_0<j_1<j_2<\ell$, if
$\eta_{p_{j_0}}\neq\eta_{p_{j_1}}$ and
$\eta_{p_{j_1}}\neq\eta_{p_{j_2}}$, then
$\nu(\eta_{p_{j_0}},\eta_{p_{j_1}})<\nu(\eta_{p_{j_1}},\eta_{p_{j_2}})$.

Let $\bar p\in{}^{e_k(\ell)}P_{n,\pi,k,\bar s}$.
Repeating the previous procedure $k$~times for subsequent
subsequences of sequences $\seq{\eta_{p_j}^i}{j<e_k(\ell)}$
for all $i<k$, we get a~one-to-one sequence $v\in{}^\ell e_k(\ell)$
such that $\set{p_{v(j)}}{j<\ell}$ fulfills the assumption
of~\autoref{K-clm} and therefore $\set{p_{v(j)}}{j<\ell}$ has an
lower bound.
For $\ell=n_i$ we get~\ref{lem:eqp} which finishes the proof
of~\autoref{Pr-ppty}.
\end{varproof}

Let $\bar n\in\incrbaire$ be such that $n_0\geq2$.
By induction define $\psi_k,h_k\in\incrbaire$ as follows:
\begin{align*}
&\psi_0(i)=i,\qquad h_k(i)=n_{\psi_k(i)},\\
&\psi_{k+1}(i)=\min\set{\ell\in\omega}
{\exists\ell'<\ell\bsp n_{\ell'}>h_k(i)\setand n_\ell\geq n_{2\ell'}\cdot h_k(i)}.
\end{align*}
Clearly, $\psi_{k+1}(i)>\psi_k(i)$ and $h_{k+1}(i)>h_k(i)$ for all $k\in\omega$.
Therefore $n_\ell>h_k(i)$ implies $\ell>k$.
Denote $\Hwf=\set{h_k}{k\in\omega}$ for all $k\in\omega$. Note that $\Lc^*_\Hwf\eqT\Nwf$.

\begin{theorem}\label{b36}
Let $\lambda$ be an uncountable regular cardinal.
Every\/ $\Pr_{\bar n}^2(\lambda)$ poset preserves
strongly $\lambda$-$\Lc^*_\Hwf$-unbounded families.
\end{theorem}

\begin{varproof}
Let $\Por$ be a~$\Pr_{\bar n^*}^2(\lambda)$ poset and assume that
$\set{x_j}{j\in J}\subseteq\baire$ is a~strongly
$\lambda$-$\Lc^*_\Hwf$-unbounded family.
Toward a contradiction, assume that there are $k_*\in\omega$ and
a~$\Por$-name $\dot\varphi$ of a member of $\Swf(\omega,h_{k_*})$
such that some $p_*\in\Por$ forces
$p_*\Vdash\lqq|\set{j\in J}{x_j\in^*\dot\varphi}|\geq\lambda\rqq$.
Then the set
$J_0=\set{j\in J}{\exists p\leq p_*\bsp(p\Vdash\lqq x_j\in^*\dot\varphi\rqq)}$
has cardinality~$\geq\lambda$.
For every $j\in J_0$, choose $p_j\leq p_*$ and $i_j\in\omega$ such that
$p_j\Vdash\lqq\forall i\geq i_j\bsp x_j(i)\in\dot\varphi(i)\rqq$.
Since $\cf(\lambda)>\omega$, we can find $i_*<\omega$ and $J_1\subseteq J_0$ of
size~$\lambda$ such that $i_j=i_*$ for all $j\in J_1$.
Since $\Por$~satisfies $\Pr^2_{\bar n^*}(\lambda)$ we can find $J_2\in[J_1]^\lambda$
and $k\geq\max\{k_*,i_*\}$ such that
$\forall i>k\bsp
\forall u\in[J_2]^{n_{i+k}}\bsp
\exists v\in[u]^{n_i}$ $\set{p_j}{j\in v}$ has a~common lower bound.
Then
\[
\tag{$\boxtimes$}\label{i>k*}
\text{$\forall i>k\bsp
\forall u\in[J_2]^{n_{2i}}\bsp
\exists v\in[u]^{n_i}$ $\set{p_j}{j\in v}$ has a~common lower bound.}
\]
Define
\[
\varphi(i)=
\begin{cases}
\set{m\in\omega}{\exists j\in J_2\bsp p_j\Vdash\lqq m\in\dot\varphi(i)\rqq},&i>k,\\
\emptyset,&i\leq k.
\end{cases}
\]

We claim that $|\varphi(i)|<h_{k+1}(i)$ for all $i>k$;
then $\varphi\in\Swf(\omega,\{h_{k+1}\})$.
Assume that for some $i>k$, $|\varphi(i)|\geq h_{k+1}(i)$ where
$h_{k+1}(i)=n_\ell\geq n_{2\ell'}\cdot h_k(i)$ and $n_{\ell'}>h_k(i)$;
then $\ell'>k$ by the remark preceding the theorem.
Since every $p_j$ forces that $|\dot\varphi(i)|\leq h_{k_*}(i)\leq h_k(i)$,
there are some $u\in[J_2]^{n_{2\ell'}}$ and a~one-to-one sequence
$\seq{m_j}{j\in u}\in{}^u\omega$ such that
$\forall j\in u\bsp p_j\Vdash\lqq m_j\in\dot\varphi(i)\rqq$.
Then by~\eqref{i>k*} there is $v\in[u]^{n_{\ell'}}$ such that an lower bound
of $\set{p_j}{j\in v}$ forces that
$|\dot\varphi(i)|\geq n_{\ell'}>h_k(i)\geq h_{k_*}(i)$.
This contradiction finishes the proof of the claim.

Since $\set{x_j}{j\in J}\subseteq\baire$ is a~strongly
$\lambda$-$\Lc^*_\Hwf$-unbounded family and $\varphi\in\Swf(\omega,\Hwf)$,
$|\set{j\in J}{x_j\in^*\varphi}|<\theta$.
Then there is $j\in J_2$ such that $x_j\not\in^*\varphi$ and then there is $i>k$
such that $x_j(i)\notin\varphi(i)$.
This leads to a~contradiction with the definition of~$\varphi(i)$ because
$i>k\geq i_*=i_j$ and $p_j\Vdash\lqq x_j(i)\in\dot\varphi(i)\rqq$.
This contradiction completes the proof of the fact that
$\set{x_j}{j\in J}$ is strongly $\lambda$-$\Lc^*_\Hwf$-unbounded in~$V^\Por$.
\end{varproof}

We are prepared to prove~\autoref{Thm:a3}.

\begin{theorem}\label{Thm:addN}

Let $\lambda$ be an uncountable regular cardinal, $|\pi|\geq\lambda$ an ordinal, and let $\Por_{\pi}=\seq{\Por_\alpha,\Qnm_\beta}{\alpha<\pi}$ be a FS iteration such that, for each $\alpha<\pi$, $\Qnm_\alpha$ is a $\Por_{\alpha}$-name of a non-trivial  $\lambda$-cc poset satisfying the property\/ $\Pr_{\bar n^*}^2(\lambda)$. Then for any cardinal $\nu\in[\lambda,|\pi|]$ with $\cf(\nu)\geq\lambda$, $\Por_\nu$ adds a strongly $\nu$-$\Lc^*_\Hwf$-unbounded family of size $\nu$ which is still strongly $\nu$-$\Lc^*_\Hwf$-unbounded in the $\Por_\pi$-extension. In particular, $\Por_\pi$ forces\/ $\add(\Nwf)\leq\lambda$ and\/ $|\pi|\leq\cof(\Nwf)$.

\end{theorem}
\begin{varproof}
The first part of the theorem is a direct consequence of~\autoref{lem:strongCohen}, \autoref{b36}, and the fact that $\Por_\pi/\Pbb_\nu$, the remaining part of the iteration from stage $\nu$, fulfills  $\Pr_{\bar n^*}^2(\lambda)$ by~\autoref{Pr^2:preservation}; and the second part is consequence of~\autoref{unbT}.
\end{varproof}

\subsection{Bonus track: More on the preservation of the additivity of the null ideal}
\

By~\autoref{lem:stsuc}, we know that $\Cor_\lambda*\Qnm$ satisfies the property \/~$\Pr_{\bar n}^2(\lambda)$ when $\Vdash_{\Cor_\lambda}\lqq\Qnm$ the property satisfies\/~$\Pr_{\bar n}^2(\lambda)\rqq$. We examine the latter to obtain a much weaker result on the preservation of $\Nwf$. 

We go over some basic facts before going into details.

\begin{fact}
Let $N\subseteq\cantor$,
$N\in\Nwf$ if and only if there is a~perfect tree $T$
such that for every $s\in T$, $\Lb([T_s])>0$ and
$N\subseteq\cantor\smallsetminus\bigcup_{m\in\omega}[T^{[m]}]$.
\end{fact}

Since Cohen forcing is forcing-equivalent to any countable atomless poset, we can use the two interpretations of Cohen forcing:
\begin{enumerate}[label=\rm(\Roman*)]
\item\label{int.I}
By Cohen forcing $\Cor$ we mean the set ${}^{<\omega}2$ ordered
by inclusion.
\item\label{int.II}

The conditions in $\Cor$ are finite trees $t\subseteq{}^{<\omega}2$
such that for some $m=m_t\in\omega$,
\begin{itemize}
\item
for every $u\in t$ there is $u'\in t\cap{}^m2$ such that
$u\subseteq u'$ and
\item $\frac{|t\cap{}^m2|}{2^m}>\frac{1}{2}$.
\end{itemize}
The ordering is defined by $t'\leq t$ iff $U_{t'}\subseteq U_t$, where
$U_t=\bigcup\set{[u]}{u\in t\cap{}^{m_t}2}$ in~$\cantor$ and
$\Lb(U_t)=\frac{|t\cap{}^{m_t}2|}{2^{m_t}}$.

\end{enumerate}
Unlike in case~\ref{int.I} the ordering in case~\ref{int.II} is not separative.
For example, trees $^{\leq n}2$, $n\in\omega$, are pairwise forcing equivalent.

\begin{theorem}
Assume that\/ $\Vdash_{\Cor_\lambda}\lqq\Qnm$ is a~forcing notion\rqq{}
and for every sequence
$\seq{(p_\xi,\dot q_\xi)}{\xi<\lambda}\in{}^\lambda(\Cor_\lambda*\Qnm)$
there are $S\in[\lambda]^\omega$, $p'_\xi\leq p_\xi$ for $\xi\in S$, and $t\in\Cor$
such that for a~sufficiently large $n\in\omega$ there are $v\in[S]^n$
and $u\in[S]^{<\omega}$ such that\/
$\bigcup_{\xi\in u\cup v}p'_\xi\in\Cor_\lambda$ and\/
$\bigcup_{\xi\in u\cup v}p'_\xi\Vdash
\lqq\set{\dot q_\alpha}{\alpha\in v}$ has a~common lower bound\rqq.
Then\/ $\Vdash_{\Cor_\lambda*\Qnm}\lqq\add(\Nwf)\leq\lambda\rqq$.
\end{theorem}

\begin{varproof}
Let $\seq{\dot T_\xi}{\xi<\lambda}$ be the canonical $\Cor_\lambda$-name
for the sequence of generic perfect subtrees of ${}^{<\omega}2$ such that
$\Vdash_{\Cor_\lambda}\lqq\Lb([\dot T_\xi])=1/2\rqq$.
Then $\dot N_\xi=\cantor\smallsetminus\bigcup_{m\in\omega}[\dot T_\xi^{[m]}]$
is a~name of a~set of measure zero.
Toward a~contradiction assume that there is $(p,\dot q)\in\Cor_\lambda*\Qnm$
such that $(p,\dot q)\Vdash\lqq\add(\Nwf)>\lambda\rqq$.
Then there exists a~name $\dot T$ for a~perfect subtree of
${}^{<\omega}2$ such that $(p,\dot q)$ forces
\begin{itemize}
\item
for every $s\in\dot T$, $\Lb([\dot T_s])>0$, and
\item
for every $\xi<\lambda$,
$[\dot T]\cap\dot N_\xi=\emptyset$,
i.e.,
$[\dot T]\subseteq\bigcup_{m\in\omega}[\dot T_\xi^{[m]}]$.
\end{itemize}
By the Baire category theorem for every $\xi<\lambda$,
\[
(p,\dot q)\Vdash\lqq\exists m\in\omega\bsp\exists s\in\dot T\bsp
\dot T_s\subseteq\dot T_\xi^{[m]}\,\rqq.
\]
Denote
\[
S_{m,s}=\set{\xi<\lambda}
{(p,\dot q)\not\Vdash\lqq\neg(\dot T_s\subseteq\dot T_\xi^{[m]}\setand
\Lb([\dot T_s])\geq2^{-m})\rqq}.
\]
Note that $\lambda=\bigcup_{m,s}S_{m,s}$.
Let $m$ and $s$ be such that $|S_{m,s}|=\lambda$.
For every $\xi\in S_{m,s}$ choose $(p_\xi,\dot q_\xi)\leq(p,\dot q)$
such that
$(p_\xi,\dot q_\xi)\Vdash\lqq\dot T_s\subseteq\dot T_\xi^{[m]}$ and
$\Lb([\dot T_s])\geq2^{-m}\,\rqq$  and $\xi\in\dom(p_\xi)$.

Let $S\in[S_{m,s}]^\omega$,
$p'_\xi\leq p_\xi$ for $\xi\in S$, and $t\in\Cor$ be such as it is promissed by the
assumptions of the lemma.
Note $U_t=\bigcup\set{[\nu]}{\nu\in t\cap{}^{m_t}2}$ is a~clopen set of measure
$\frac{|t\cap{}^{m_t}2|}{2^{m_t}}>\frac{1}{2}$ and
we can assume that $m\geq m_t$.
Let $\ell>1$ be such that
\[\frac{|t\cap{}^{m_t}|}{2^{m_t}}>\frac12+2^{-\ell}.\]
For $n=2^\ell$ find
$v\in[S]^n$ and $u\in[S]^{<\omega}$ such that
$p'=\bigcup_{\xi\in u\cup v}p'_\xi$ is a~condition and $p'\Vdash
\lqq\set{\dot q_\alpha}{\alpha\in v}$ has a~common lower bound\rqq\ and 
let $\dot q'$ be the $\Cor_\lambda$-name for the the common
extension of $\set{\dot q_\alpha}{\alpha\in v}$.

Choose enumerations ${}^\ell2=\set{\nu_j}{j<2^\ell}$ and $v=\set{\xi_j}{j<2^\ell}$.
For every $j<2^\ell$ let $t_j\in\Cor$ be such that $t_j\leq t$, $m_{t_j}=m+\ell$, and
$U_{t_j}=U_t\smallsetminus\bigcup_{\eta\in{}^m2}[\eta^\frown\nu_j]$.
This is possible because
\[\Lb(U_{t_j})=\Lb(U_t)\cdot(1-2^{-\ell})\geq\Lb(U_t)-2^{-\ell}>1/2.\]
Clearly,
$U_{t_j}^{[m]}\subseteq\cantor\smallsetminus\bigcup_{\eta\in{}^m2}[\eta^\frown\nu_j]$
and $\bigcap_{j<2^\ell}U_{t_j}^{[m]}=\emptyset$.
We define an extension $p''$ of~$p'$ as follows:
$\dom(p'')=\dom(p')$, $p''(\xi)=p'(\xi)$ for $\xi\notin v_0$, and
$p''(\xi_j)=t_j$ for $j<2^\ell$.
Now, $(p'',\dot q')\leq(p',\dot q')\leq(p,\dot q)$ and
$(p',\dot q')\Vdash\lqq\Lb\big([\dot T_s]\big)\geq2^{-m}\setand
[\dot T_s]\subseteq\bigcap_{\xi\in v}[\dot T_\xi^{[m]}]\subseteq
\bigcap_{j<2^\ell}U_{t_j}^{[m]}=\emptyset\rqq$.
This is a~contradiction.
\end{varproof}

\section{Applications I: Separating the right side of Cicho\'n's diagram along with \texorpdfstring{$\efrak_2^\const$}{}}\label{sec:s4}

The purpose of this section is to provide applications of our results obtained in the prior section when forcing iterations. In particular, we
perform forcing constructions by using the method of matrix iterations with ultrafilters
from~\cite{BCM} to prove~\eqref{Seplef:1}--\eqref{Seplef:2} where

We start with proof of~\eqref{Seplef:1}:

\begin{theorem}\label{e0}
Let $\lambda_1\leq\lambda_2\leq\lambda_3\leq\lambda_4\leq\lambda_5$ be uncountable regular cardinals, and $\lambda_6$~a~cardinal such that $\lambda_6\geq\lambda_5$ and\/ $\cof([\lambda_6]^{<\lambda_i})=\lambda_6 = \lambda_6^{\aleph_0}$ for $1\leq i\leq 4$. Then we can built a ccc poset forcing
\begin{enumerate}[label=\rm (\alph*)]
\item\label{e1:a}
$\cfrak = \lambda_6$;
\item\label{e1:c}
$\Nwf\eqT\Cv_{[\lambda_6]^{<{\lambda_1}}}$, $\Esf_2\eqT\Cv_{[\lambda_6]^{<{\lambda_2}}}$, $\Cv_{\Nwf}^{\perp}\eqT\Cv_{[\lambda_6]^{<{\lambda_3}}}$, and $\baire\eqT\Cv_{[\lambda_6]^{<{\lambda_4}}}$;
\item\label{e1:d}
$\lambda_5\leqT\Cv_\Mwf$, $\lambda_6\leqT \Cv_\Mwf$, and $\Ed\leqT \lambda_6\times\lambda_5$.
\end{enumerate}
In particular, it forces
\begin{multline*}
\add(\Nwf)=\lambda_1\leq\efrak_K^\const=\efrak_2^\const=\lambda_2\leq\cov(\Nwf)=\lambda_3\leq\bfrak=\lambda_4\leq\\
\leq\non(\Mwf)=\lambda_5\leq\cov(\Mwf)=\lambda_6=\cfrak.
\end{multline*}
\end{theorem}

\begin{varproof}
For each $\rho<\lambda_6\lambda_5$ denote $\eta_\rho:=\lambda_6\rho$. Fix a bijection $g=(g_0, g_1,g_2):\lambda_6\to\{0,1,2,3\}\times\lambda_6\times\lambda_6$ and fix a function $t\colon\lambda_6\lambda_5\to\lambda_6$ such that, for any $\alpha<\lambda_6$, the set $\set{\rho<\lambda_6\lambda_5}{t(\rho)=\alpha}$ is cofinal in $\lambda_6\lambda_5$.

We below construct the ${<}\lambda_4$-uf-extendable matrix iteration $\Por_{\gamma,\pi}$ with $\gamma=\lambda_6$ and $\pi=\lambda_6\lambda_6\lambda_5$. First set,
\begin{enumerate}[label=\rm (C\arabic*)]
\item\label{(C1)}
$\Delta(\alpha):=\alpha+1$ and $\Qnm^*_{\alpha}=\Cor_{\alpha}$ for $\alpha\leq\lambda_5$.
\end{enumerate}
Let us define the matrix iteration at each $\xi=\lambda_\rho+\varepsilon$ for $\rho<\lambda_6\lambda_5$ and $\varepsilon<\lambda_6$ as follows. Denote\footnote{We think of $X_1$ as the set of Borel codes of Borel sets with measure zero.}
\begin{align*}
\Qor^*_0 & := \Loc, & \Qor^*_1 &:= \Por^2, & \Qor^*_2 & := \Bor, & \Qor^*_3 &:=\Dor,\\
X_0 & := \baire, & X_1 & := {}^{\omega}2, & X_2 &:= \Bwf(\cantor)\cap \Nwf, & X_3 & := \baire.
\end{align*}
For $j<4$, $0<\rho<\lambda_6\lambda_5$ and $\alpha<\lambda_5$, choose
\begin{enumerate}[label=(E$j$)]
\item\label{Ej}
a collection $\set{\Qnm_{j,\alpha,\zeta}^\rho}{\zeta<\lambda_6}$ of nice $\Por_{\alpha,\eta_\rho}$-names for posets of the form $(\Qor^*_j)^N$ for some transitive model $N$ of $\thzfc$ with $|N|<\lambda_{j+1}$
such that, for any $\Por_{\alpha,\eta_\rho}$-name $\dot F$ of a subset of $X_j$ of size ${<}\lambda_{j+1}$, there is some $\zeta<\lambda_6$ such that, in $V_{\alpha,\eta_\rho}$, $\Qnm^\rho_{j,\alpha,\zeta} = (\Qor^*_j)^N$ for some $N$ containing $\dot F$ (we explain later why this is possible),
\end{enumerate}
and set
\begin{enumerate}[label=\rm (C\arabic*)]
\setcounter{enumi}{1}
\item
$\Delta(\xi):=t(\rho)$ and $\Qnm_{\xi}:=\Eor^{V_{\Delta(\xi),\xi}}$ when $\xi=\eta_\rho$;
\item
$\Delta(\xi):=g_1(\varepsilon)$ and $\Qnm_{\xi}:=\Qnm^{\rho}_{g(\varepsilon)}$ when $\xi=\eta_\rho+1+\varepsilon$ for some $\varepsilon< \lambda_6$.
\end{enumerate}
According to~\autoref{f4}, the above settles the construction of $\Por$ as an ${<}\lambda_4$-uf-extendable matrix iteration by \autoref{exm:ufl} and \autoref{Pr_K-uf}. First, observe that $\Por$ is ccc. It is also clear that $\Por$ forces~\ref{e1:a} by assumption $\lambda_5 = \lambda_5^{\aleph_0}$. We now prove that $\Por$ forces what we want:

$\Por$ forces~\ref{e1:c}. Why? Indeed:
\begin{enumerate}[label=\rm(\faPagelines$_\arabic*$)]
\item\label{c:0}
$\Por$ forces $\Esf_2\eqT\Cv_{[\lambda_6]^{<{\lambda_2}}}$. Firstly, $\Cv_{[\lambda_6]^{<{\lambda_2}}}\leqT\Esf_2^k$ is forced by~\autoref{b4} because, for each $\xi<\pi$, $\Por_{\gamma,\xi}$ forces that $\Qnm_{\gamma,\xi}$ is $\lambda_2$-$\Esf_2^k$-good. In fact:
\begin{itemize}
\item The cases $\xi<\lambda_6$ and $\xi=\eta_\rho$ for $\rho>0$ follow by~\autoref{c11};
\item when $\xi=\eta_\rho+1+\varepsilon$ for some $\rho>0$ and $\varepsilon<\lambda_6$, we split into four subcases:
\begin{itemize}
\item the case $g_0(\varepsilon)=0$ is clear by~\autoref{b22};
\item when $g_0(\varepsilon)=1$ it follows by~\autoref{b22};
\item when $g_0(\varepsilon)=2$, it follows by~\autoref{c10} (recall $\Bor$ is $\sigma$-$n$-linked); and
\item when $g_0(\varepsilon)=3$, use~\autoref{c11}.
\end{itemize}
\end{itemize}
Since $\Por$ forces $\Cv_{[\lambda_6]^{<{\lambda_2}}}\leqT\Esf_2^k$ and $\thzfc$ proves $\Esf_2^k\leqT\Esf_2$ by~\autoref{cxn:e_2}, we obtain $\Por$ forces $\Cv_{[\lambda_6]^{<{\lambda_2}}}\leqT\Esf_2$. On the other hand, we show that $\Por$ forces $\Esf_2\leqT\Cv_{[\lambda_6]^{<{\lambda_2}}}$. Let $\dot A$ be a $\Por$-name for a subset of $\cantor$ of size ${<}\lambda_2$. By employing~\autoref{realint} we can can find $\alpha<\lambda_6$ and $\rho<\lambda_6\lambda_5$ such
that $\dot A$ is $\Por_{\alpha,\eta_\rho}$-name. By~(E0), we can find a $\zeta<\lambda_6$ and a $\Por_{\alpha,\eta_\rho}$-name $\dot N$ of a transitive model
of $\thzfc$ of size ${<}\lambda_2$ such that $\Por_{\alpha,\eta_\rho}$ forces that $\dot N$ contains $\dot A$ as a subset and $(\Por^2)^{\dot N}=\Qnm_{1,\alpha,\zeta}^\rho$, so the
generic predictor added by $\Qnm_{\xi}=\Qnm_{g(\varepsilon)}^\rho$ predicts all the reals in $\dot A$ where $\varepsilon:=g^{-1}(1,\alpha,\zeta)$ and $\xi=\eta_\rho+1+\varepsilon$. Then, by applying~\autoref{b12}, $\Por$ forces $\Esf_2\leqT\Cv_{[\lambda_6]^{<{\lambda_2}}}$. 
\end{enumerate}

Find a function $h\in\baire$ that converges to infinity. Then by using~\autoref{Pr-ppty}~\ref{Pr-ppty:b}, we can find a $\bar\rho$ such that $\Por^2$ is $\sigma$-$\bar\rho$-linked. Hence, by applying ~\autoref{goodrho}, there is $\Hwf_{\bar\rho,h}=\set{h_n}{n\in\omega}\subseteq\baire$ such that $\Por^2$ is $\Lc^*_{\Hwf_{\bar\rho,h}}$-good. 

\begin{enumerate}[start=2,label=\rm(\faPagelines$_\arabic*$)]
\item $\Por$ forces $\Cv_{[\lambda_6]^{<{\lambda_1}}}\eqT\Lc^*_{\Hwf_{\bar\rho,h}}$: Firstly, $\Por$ forces $\Cv_{[\lambda_6]^{<{\lambda_1}}}\leqT\Lc^*_{\Hwf_{\bar\rho,h}}$ since, it can be checked as in~\ref{c:0} that all iterands are $\lambda_1$-$\Lc^*_{\Hwf_{\bar\rho,h}}$-good (see~\autoref{rem:rancent}) so, by~\autoref{b4} $\Por$ forces $\Cv_{[\lambda_6]^{<{\lambda_1}}}\leqT\Lc^*_{\Hwf_{\bar\rho,h}}$; and secondly $\Lc^*_{\Hwf_{\bar\rho,h}}\leqT\Cv_{[\lambda_6]^{<{\lambda_4}}}$ is forced as in~\ref{c:0}. 

\item $\Por$ forces $\Cv_{[\lambda_6]^{<{\lambda_4}}}\eqT\baire$: $\Por$ forces $\Cv_{[\lambda_6]^{<{\lambda_4}}}\leqT\baire$ by~\autoref{mainpres} because the matrix iteration is ${<}\lambda_4$-uf-extendable and $\baire\leqT\Cv_{[\lambda_6]^{<{\lambda_4}}}$ is forced as in~\ref{c:0}.

\item $\Por$ forces $\Cv_{\Nwf}^{\perp}\eqT\Cv_{[\lambda_6]^{<{\lambda_3}}}$. This is basically the same argument as in~\ref{c:0}.

\end{enumerate}
$\Por$ forces~\ref{e1:d}:
Given that $\cf(\pi)=\lambda_5$, $\Por$ forces that $\lambda_5\leqT\Cv_\Mwf$ by applying~\autoref{lem:strongCohen}, and, by~\autoref{matsizebd}, $\Por$ forces $\lambda_6\leqT\Cv_\Mwf$. It is yet to be proved that $\Por$ forces that $\Ed\leqT \lambda_6\times\lambda_6\lambda_5$ (because $\lambda_6\lambda_5\eqT \lambda_5$). For this purpose, for each $\rho<\lambda_6\lambda_5$ denote by $\dot e_\rho$ the $\Por_{\Delta(\eta_\rho),\eta_\rho+1}$-name of the eventually different real over $V_{t(\rho),\eta_\rho}$ added by $\Qnm_{t(\rho),\eta_\rho}$. In $V_{\gamma, \pi}$, we are going to define maps $\Psi_-:\baire\to\lambda_6\times\lambda_6\lambda_5$ and $\Psi_+:\lambda_6\times\lambda_6\lambda_5\to\baire$ such that, for any $x\in\baire$ and for any $(\alpha,\rho)\in\lambda_6\times\lambda_6\lambda_5$, if $\Psi_-(x)\leq(\alpha,\rho)$ then $x\neq^\infty\Psi_+(\alpha,\rho)$.

For $x\in V_{\lambda_6,\pi}\cap\baire$, we can find $\alpha_x<\lambda_6$ and $\rho_x<\lambda_6\lambda_5$ such that $x\in V_{\alpha_x,\eta_{\rho_x}}$, so put $\Psi_-(x):=(\alpha_x,\rho_x)$; for $(\alpha, \rho)\in\lambda_6\times\lambda_6\lambda_5$, find some $\rho'<\lambda_6\lambda_5$ such that $\rho'\geq\rho$ and $t(\rho')=\alpha$, and define $\Psi_+(\alpha,\rho):=\dot e_{\rho'}$. It is clear that $(\Psi_-,\Psi_+)$ is
the required Tukey connection.
\end{varproof}

We are ready to proceed with the next application, that is, with proof of~\eqref{Seplef:2}.

\begin{theorem}\label{c1}
Let $\lambda_1\leq\lambda_2\leq\lambda_3\leq\lambda_4\leq\lambda_5$ be uncountable regular cardinals, and $\lambda_6$~a~cardinal such that $\lambda_6\geq\lambda_5$ and $\cof([\lambda_6]^{<\lambda_i})=\lambda_6 = \lambda_6^{\aleph_0}$ for $1\leq i\leq 4$. Then we can built a ccc poset forcing
\begin{enumerate}[label=\rm (\alph*)]
\item\label{c1:a}
$\cfrak = \lambda_6$;
\item\label{c1:c}
$\Esf_\omega\eqT\Esf_2\eqT\Cv_{[\lambda_6]^{<{\lambda_1}}}$, $\Nwf\eqT\Cv_{[\lambda_6]^{<{\lambda_2}}}$, $\Cv_{\Nwf}^{\perp}\eqT\Cv_{[\lambda_6]^{<{\lambda_3}}}$, and\/ $\baire\eqT\Cv_{[\lambda_6]^{<{\lambda_4}}}$;
\item\label{c1:d}
$\lambda_5\leqT\Cv_\Mwf$, $\lambda_6\leqT \Cv_\Mwf$, and\/ $\Ed\leqT \lambda_6\times\lambda_5$.
\end{enumerate}
In particular, it forces
\begin{multline*}
\efrak^\const=\efrak_K^\const=\efrak_2^\const=\lambda_1\leq\add(\Nwf)=\lambda_2\leq\cov(\Nwf)=\lambda_3\leq\bfrak=\lambda_4\leq\\
\leq\non(\Mwf)=\lambda_5\leq\cov(\Mwf)=\lambda_6=\cfrak.
\end{multline*}
\end{theorem}

\begin{varproof}
We proceed as in the proof of \autoref{e0}. Let $\eta_\rho$ ($\rho<\lambda_6\lambda_5$) and $t$ be as there, and fix a bijection $g=(g_0, g_1,g_2):\lambda_6\to\{0,1,2,3\}\times\lambda_6\times\lambda_6$. We will construct a ${<}\lambda_4$-uf-extendable matrix iteration $\Por = \Por_{\gamma,\pi}$ with $\gamma = \lambda_6$ and $\pi = \lambda_6\lambda_6\lambda_5$. 

First let \ref{(C1)} be as in the proof~\autoref{e0} and define the matrix iteration at each $\xi=\eta_\rho+\varepsilon$ for $0<\rho<\lambda_6\lambda_5$ and $\varepsilon<\lambda_6$ as follows. Denote
\begin{align*}
\Qor^*_0 & := \Por^\omega, & \Qor^*_1 &:= \Loc, & \Qor^*_2 & := \Bor, & \Qor^*_3 &:=\Dor,\\
X_0 & := \baire, & X_1 & := \baire, & X_2 &:= \Bwf(\cantor)\cap \Nwf, & X_3 & := \baire.
\end{align*}
For $j<4$, $0<\rho<\lambda_6\lambda_5$ and $\alpha<\lambda_5$, choose
\begin{enumerate}[label=(E$j$)]
\item
a collection $\set{\Qnm_{j,\alpha,\zeta}^\rho}{\zeta<\lambda_6}$ of nice $\Por_{\alpha,\eta_\rho}$-names for posets of the form $(\Qor^*_j)^N$ for some transitive model $N$ of $\thzfc$ with $|N|<\lambda_{j+1}$
such that, for any $\Por_{\alpha,\eta_\rho}$-name $\dot F$ of a subset of $X_j$ of size ${<}\lambda_{j+1}$, there is some $\zeta<\lambda_6$ such that, in $V_{\alpha,\eta_\rho}$, $\Qnm^\rho_{j,\alpha,\zeta} = (\Qor^*_j)^N$ for some $N$ containing $\dot F$ (we explain later why this is possible),
\end{enumerate}
and set
\begin{enumerate}[label=\rm (C\arabic*)]
\setcounter{enumi}{1}
\item
$\Delta(\xi):=t(\rho)$ and $\Qnm_{\xi}:=\Eor^{V_{\Delta(\xi),\xi}}$ when $\xi=\eta_\rho$;
\item
$\Delta(\xi):=g_1(\varepsilon)$ and $\Qnm_{\xi}:=\Qnm^{\rho}_{g(\varepsilon)}$ when $\xi=\eta_\rho+1+\varepsilon$ for some $\varepsilon<\lambda_6$.
\end{enumerate}
The construction is indeed a ${<}\lambda_4$-uf-extendable iteration. We just prove that $\Por$ forces $\Rbf_\omega\eqT\Esf_2\eqT\Cv_{[\lambda_6]^{<{\lambda_1}}}$, as the rest can be proved as in \autoref{e0}. To see that $\Por$ forces $\Cv_{[\lambda_6]^{<{\lambda_1}}}\leqT\Esf_2^k$, it suffices to check that for each $\xi<\pi$, $\Por_{\gamma,\xi}$ forces that $\Qnm_{\gamma,\xi}$ is $\lambda_2$-$\Esf_2^k$-good. In fact:
\begin{itemize}
\item The cases $\xi<\lambda_6$ and $\xi=\eta_\rho$ for $\rho>0$ follow by~\autoref{c11};
\item when $\xi=\eta_\rho+1+\varepsilon$ for some $\rho>0$ and $\varepsilon<\lambda_6$, we split into four subcases:
\begin{itemize}
\item the case $g_0(\varepsilon)=0$ is clear by~\autoref{b22};
\item when $g_0(\varepsilon)=1$ it follows by~\autoref{c10} (recall $\Loc$ is $\sigma$-$n$-linked);
\item when $g_0(\varepsilon)=2$, it follows by~\autoref{c10} (recall $\Bor$ is $\sigma$-$n$-linked); and
\item when $g_0(\varepsilon)=3$, use~\autoref{c11}.
\end{itemize}
\end{itemize}
Hence, by~\autoref{b4} we obtain $\Por$ forces $\Cv_{[\lambda_6]^{<{\lambda_1}}}\leqT\Esf_2^k$ and since $\thzfc$ proves $\Esf_2^k\leqT\Esf_2$ by using~\autoref{cxn:e_2}, so it is forced $\Cv_{[\lambda_6]^{<{\lambda_1}}}\leqT\Esf_2$. On the other hand, $\Por$ forces $\Rbf_\omega\leqT\Cv_{[\lambda_6]^{<{\lambda_1}}}$ is proved as in~\ref{c:0} in the proof of~\autoref{e0}. Thereby, $\Por$ forces $\Cv_{[\lambda_6]^{<{\lambda_1}}}\leqT\Esf_2\leqT\Rbf_\omega\leqT\Cv_{[\lambda_6]^{<{\lambda_1}}}$ by~\autoref{t2}.
\end{varproof}

The construction of~\cite[Thm.~2.43]{KST} for the alternative order of the left side of Cicho\'n’s diagram can also be adapted to the context of the earlier theorems.  By simply switching the order of the values of $\bfrak$ and $\cov(\Nwf)$, we are forcing $\bfrak= \lambda_3\leq \cov(\Nwf) =\lambda_4$ rather than $\cov(\Nwf) = \lambda_3\leq\bfrak= \lambda_4$. This construction used the method of iterations with finitely additive measures (fams), which is more powerful than the method preceding ultrafilter-limits. The latter was introduced by~Shelah~\cite{ShCov} to prove the consistency of $\thzfc$ with $\cf(\cov(\Nwf)) = \omega$. Recently, these works were formalized the general framework of iterations with fam-limits in~\cite{CMU}.

We below use the forcing construction done by Mej\'ia~\cite[Thm.~4.1]{modKST} which weakening of the hypothesis $\thgch$ originally used in~\cite[Thm.~2.43]{KST} to establish~\eqref{Seplef:3}. 

\begin{theorem}\label{c2}
Let $\lambda_1\leq\lambda_2\leq\lambda_3\leq\lambda_4\leq\lambda_5$ be uncountable regular cardinals, and let $\lambda_6 = \lambda_6^{\aleph_0} \geq \lambda_5$ be a cardinal such that\/ $\cof([\lambda_6]^{<\lambda_i}) = \lambda_6$ for\/ $1 \leq i \leq 3$.
Further assume that one of the following holds:
\begin{enumerate}[label=\rm(\alph*)]
\item
$\lambda_3=\lambda_4$.
\item
both $\lambda_4$ and $\lambda_5$ are\/ $\aleph_1$-inaccessible,
and there is some cardinal $\lambda <\lambda_4$ such that\/ $\forall \alpha<\lambda_3\colon |\alpha|\leq\lambda $, and whenever\/ $2^\lambda <\lambda_6$, $\lambda ^{<\lambda } = \lambda $.
\end{enumerate}
Then, we can construct a FS iteration of length $\lambda_6$ of ccc posets forcing
\begin{multline*}
\add(\Nwf)=\lambda_1\leq\efrak_2^\const=\lambda_2\leq\bfrak=\lambda_3\leq\cov(\Nwf)=\lambda_4\leq\\
\leq\non(\Mwf)=\lambda_5\leq\cov(\Mwf)=\cfrak=\lambda_6.
\end{multline*}
\end{theorem}

\begin{varproof}
We force with the poset from~\cite[Thm.~4.1]{modKST}, which we call $\Qor$. We review without going into the details of how $\Qor$ is constructed for completenes's sake. This forcing construction uses a FS iteration $\seq{\Por_\alpha,\Qnm_\alpha}{\alpha<\pi}$ of length $\pi:=\lambda_6+\lambda_6$ with fam-limits,
going through of the following ccc posets:
\begin{enumerate}[label=\rm(P\arabic*)]
\item
restrictions of $\tilde\Eor$ (a variant of $\Eor$ to add an eventually different real) to $V^{\Por'_{\alpha}}$ for some $\Por'_\alpha \subsetdot \Por_\alpha$ of size ${<}\lambda_5$; 
\item
restrictions of $\Bor$ to $V^{\Por'_{\alpha}}$ for some $\Por'_\alpha \subsetdot \Por_\alpha$ of size ${<}\lambda_4$;
\end{enumerate}
and use small transitive models as in~\autoref{e0} for
\begin{enumerate}[label=\rm(P\arabic*)]
\setcounter{enumi}{2}
\item
all $\sigma$-centered subposets of $\Dor$ of size~${<}\lambda_3$.
\item\label{smallP^2} 
all $\sigma$-linked posets of $\Por^2$ of size~${<}\lambda_2$;
\item
all $\sigma$-linked posets of $\Loc$ of size~${<}\lambda_1$;
\end{enumerate}
As in the proof of~\autoref{e0}, find a function $h\in\baire$ that converges to infinity, so by using~\autoref{Pr-ppty}~\ref{Pr-ppty:b}, there is a $\bar\rho$ such that $\Por^2$ is $\sigma$-$\bar\rho$-linked. Therefore, by applying ~\autoref{goodrho}, there is $\Hwf_{\bar\rho,h}=\set{h_n}{n\in\omega}\subseteq\baire$ such that $\Por^2$ is $\Lc^*_{\Hwf_{\bar\rho,h}}$-good. 

We prove that $\Qor$ is as required.
We just are left with seeing $\Qor$ forces $\efrak_2^\const=\lambda_2$ and $\add(\Nwf)=\lambda_1$ since the rest is exactly as in the proof of~\cite[Thm.~4.1]{modKST}. On the one hand, $\Qor$ forces $\efrak_2^\const\leq\lambda_2$ and $\add(\Nwf)\leq\lambda_1$ by~\autoref{b4} because all iterands are $\Esf_2^k$-good and $\Lc^*_{\Hwf_{\bar\rho,h}}$-good; on the other hand, $\efrak_2^\const\geq\lambda_2$ and $\add(\Nwf)\geq\lambda_1$ is forced as~\autoref{e0}. 
\end{varproof}

Just as in~\autoref{c2}, we can also prove the following.

\begin{theorem}
Let $\lambda_1\leq\lambda_2\leq\lambda_3\leq\lambda_4\leq\lambda_5$ be uncountable regular cardinals, and let $\lambda_6 = \lambda_6^{\aleph_0} \geq \lambda_5$ be a cardinal such that\/ $\cof([\lambda_6]^{<\lambda_i}) = \lambda_6$ for $1 \leq i \leq 3$.
Further assume that one of the following holds:
\begin{enumerate}[label=\rm(\alph*)]
\item $\lambda_3=\lambda_4$.
\item both $\lambda_4$ and $\lambda_5$ are\/ $\aleph_1$-inaccessible,
and there is some cardinal $\lambda <\lambda_4$ such that\/ $\forall \alpha<\lambda_3\colon |\alpha|\leq\lambda $, and whenever\/ $2^\lambda <\lambda_6$, $\lambda ^{<\lambda } = \lambda $.
\end{enumerate}
Then, we can construct a FS iteration of length $\lambda_6$ of ccc posets forcing
\begin{multline*}
\efrak^\const=\efrak_2^\const=\lambda_1\leq\add(\Nwf)=\lambda_2\leq\bfrak=\lambda_3\leq\cov(\Nwf)=\lambda_4\leq\\
\leq\non(\Mwf)=\lambda_5\leq\cov(\Mwf)=\cfrak=\lambda_6.
\end{multline*}
\end{theorem}

\section{Applications II: Expanding Cicho\'n's maximum}\label{sec:subm}

Using the intersections with elementary submodels approach from~\cite{GKMS}, we aim to display~\autoref{Thm:a0}-\ref{Thm:a2}, that is, we force Cichon’s maximum simultaneously along with  $\mathfrak{e}_2^{\mathrm{cons}}$ and $\mathfrak{v}_2^{\mathrm{cons}}$ with distinct values. We adhere to the presentation from~\cite[Sec.~4 \&~5]{CM22}, which is close to work. To this end, set a regular cardinal $\chi$ that is sufficiently large. The objective is to find some $N\preceq H_\chi$, closed under countable sequences, such that $\Por\cap N$ is the desired ccc poset, where $\Por$ is the ccc poset constructed in \autoref{e0} for large enough $\lambda_i$.

Taking into account that the method for determining $N$ follows the principles described in \cite[Sec.~3]{GKMS} and \cite[Sec.~5]{CM22}. In this work, it is clear that an explicit construction of the model $N$ is unnecessary, as it leads to something repetitive in this entire process in order to show~\autoref{Thm:a0}-\ref{Thm:a2}. Instead, we will discuss the impact of $\Por\cap N$ on the Tukey connections induced by $\Por$, which holds greater importance than the model's construction. Specifically, when $\Rbf$ constitutes a definable relational system of the reals (\autoref{def:defrel}) or possesses sufficient absoluteness with parameters in $N$, and if $K = \la A,B,\leqtr\ra\in N$ represents a relational system, then \[\Vdash_\Por \Rbf\leqT K\Rightarrow\ \Vdash_{\Por\cap N} \Rbf\leqT K\cap N\ \text{(by~\autoref{KcapN})}.\] Here, $K\cap N:=\la A\cap N, B\cap N,\leqtr\ra$, and both $K$ and $K\cap N$ are considered fixed within the ground model; that is, they are not subject to reinterpretation in generic extensions. A corresponding result applies for the case where  $K\leqT \Rbf$. Consequently, the values that $\Por\cap N$ forces to $\bfrak(\Rbf)$ and $\dfrak(\Rbf)$ are determined by the values of $\bfrak(K\cap N)$ and $\dfrak(K\cap N)$, which do not change in ccc forcing extensions for the $K$ of our interest. In~\cite[Sec.~4~\&~5]{CM22}, a comprehensive explanation is provided regarding how such an $N$ should be constructed to yield the desired values of $\bfrak(K\cap N)$ and $\dfrak(K\cap N)$.

In general terms, let $\kappa$ be an uncountable regular cardinal, $\chi$ a large enough regular cardinal, $N\preceq H_\chi$ a \emph{${<}\kappa$ closed model}, i.e.\ $N^{<\kappa}\subseteq N$, and let $\Por\in N$ be a $\kappa$-cc poset. Then: 

\begin{lemma}[{see e.g~\cite[Lem.~59]{Brmodtec}}]  \label{submodelbasic}
\ 
\begin{enumerate}[label=\rm(\arabic*)]
\item For every antichain $A \subseteq \Por$, $A \in N$ if and only if $A \subseteq N$. 
\item $\Por \cap N$ is $\kappa$-cc.
\item $\Por \cap N  \subsetdot\Por$.
\end{enumerate}
\end{lemma}
Even more, there is a one-to-one correspondence between the (nice) $\Por\cap N$-names of members of $H_\kappa$, and the $\Por$-names $\tau\in N$ of members of $H_\kappa$, in particular, we have the same correspondence between the nice names of reals. Moreover, if $G$ is $\Por$-generic over $V$, then $H_\kappa^{N[G]} = H_\kappa^{V[G]}\cap N[G]$ (see e.g~\cite[Fact 4.3]{CM22}). We also have absoluteness results: 

\begin{lemma}
If $p\in\Por\cap N$, $\varphi(\bar x)$ is a sufficiently absolute formula and $\bar \tau\in N$ is a finite sequence of $\Por$-names of members of $H_\kappa$, then
\[p\Vdash_\Por\varphi(\bar\tau) \text{ iff }p\Vdash_{\Por\cap N}\varphi(\bar\tau).\]    
\end{lemma}

Let $\Rbf=\la X,Y,\sqsubset\ra$ and $K=\la A,B,\leqtr\ra$ be relational systems. Note that
\begin{align*}
    K\leqT \Rbf \text{ iff } & \text{there is a sequence $\la x_a:\, a\in A\ra$ in $X$ such that}\\
     & \forall\, y\in Y\ \exists\, b_y\in B\ \forall\, a\in A\colon a \nleqtr b_y \imp x_a \nsqsubset y;\\[1ex]
    \Rbf \leqT K \text{ iff } & \text{there is a sequence $\la y_b:\, b\in B\ra$ in $Y$ such that}\\
    & \forall\, x\in X\ \exists\, a_x\in A\ \forall\, b\in B\colon a_x \leqtr b \imp x\sqsubset y_b.
\end{align*}

In the first case, we say that the sequence $\seq{x_a}{a\in A}$ \emph{witnesses $K\leqT\Rbf$}, while in the second case we say that the sequence $\seq{y_b}{b\in B}$ \emph{witnesses $\Rbf\leqT K$}.

We focus on the case in which $K$ represents a fixed relational system within the ground model (meaning that it is interpreted as the same relational system across any model of $\thzfc$). In contrast, $\Rbf$ may vary depending on its interpretation within different models, particularly when $\Rbf$ is a gPrs. When considering a poset $\Por$, it is evident that $\Vdash_\Por K\leqT \Rbf$ if and only if there is a sequence $\seq{\dot x_a}{ a\in A}$ of members of $X$ such that $\Por$ forces that $\seq{ \dot x_a}{ a\in A}$ witnesses $K\leqT\Rbf$. In this context, we denote that the sequence $\seq{\dot x_a}{ a\in A}$ \emph{witnesses $\Vdash_\Por K\leqT \Rbf$}. Similarly, one can establish the definition for “$\seq{\dot y_b}{b\in B}$ \emph{witnesses $\Vdash_\Por\Rbf\leqT K$}.”

In this situation, it is known the following:

\begin{lemma}[{\cite[Lem.~4.5]{CM22}}]\label{KcapN}
Let $\kappa$ be an uncountable regular cardinal, $\chi$ a large enough regular cardinal, and let $N\preceq H_\chi$ be a ${<}\kappa$-closed model. Assume that\/ $\Por\in N$ is a $\kappa$-cc poset, $K=\la A,B,{\leqtr}\ra\in N$ is a relational system, and\/ $\Rbf$ is a definable relational system of the reals with parameters in $N$.
Then:
\begin{enumerate}[label = \rm (\alph*)]
\item
If\/ $\seq{\dot y_b}{b\in B}\in N$ witnesses\/ $\Vdash_{\Por}\lqq \Rbf\leqT K\rqq$ then\/ $\seq{\dot y_b}{b\in B\cap N}$ witnesses\/ $\Vdash_{\Por\cap N}\lqq \Rbf \leqT K\cap N\rqq$.
\item
If\/ $\seq{\dot x_a}{a\in A}\in N$ witnesses\/ $\Vdash_{\Por}\lqq K\leqT \Rbf\rqq$ then $\seq{\dot x_a}{a\in A\cap N}$ witnesses\/ $\Vdash_{\Por\cap N}\lqq K\cap N\leqT\Rbf\rqq$.
\end{enumerate}
\end{lemma}

We are finally ready to proceed with the proof of~\autoref{Thm:a0}. Denote the relational systems (some introduced in~\autoref{b16} and~\ref{b30}) 
\[\Rbf_i:=\left\{\begin{array}{ll}
        \Lc^*  & \text{if $i=1$;}\\
         \Esf_2  & \text{if $i=2$;} \\
         \Cv_{\Nwf}^{\perp} & \text{if $i=3$;}\\
        \baire  & \text{if $i=4$;}\\
    \end{array}\right.\]
and $\Rbf_2^*=\Esf_2^k$. Recall that $\thzfc$ proves $\Esf_2^k\leqT\Esf_2$ and $\Cv_\Mwf \leqT \Ed$.
\begin{varproof}[Proof of~\autoref{Thm:a0}]
Without loss of generality, we assume that GCH holds above some regular cardinal $\lambda^-_1>\lambda$.
Let $\lambda^-_i$ and $\lambda_i$ be ordered as in~\autoref{fig:cichoncoll}, where the inequalities between them are strict. Next let $\Por$ be a ccc poset obtained from~\autoref{e0}. According to the proof, $\Por$ forces:

\begin{figure}[ht]
\centering
\begin{tikzpicture}[xscale=2/1]
\footnotesize{
\node (addn) at (0,3){$\lambda_1$};
\node (covn) at (0,6){$\lambda_3$};
\node (e2) at (0.7,4.5){$\lambda_2$};
\node (nonn) at (3.8,3) {$\lambda_6$} ;
\node (nonsn) at (3.2,3) {$\lambda_6$} ;
\node (cfn) at (3.8,6) {$\lambda_6$} ;
\node (addm) at (1.4,3) {$\bullet$} ;
\node (covm) at (2.6,3) {$\lambda_6$} ;
\node (nonm) at (1.4,6) {$\lambda_5$} ;
\node (cfm) at (2.6,6) {$\bullet$} ;
\node (b) at (1.4,4.5) {$\lambda_4$};
\node (d) at (2.6,4.5) {$\lambda_6$};
\node (c) at (5,6) {$\lambda_6$};

\draw[gray]
      (covn) edge [->] (nonm)
      (nonm)edge [->] (cfm)
      (cfm)edge [->] (cfn)
      (addn) edge [->]  (addm)
      (addn) edge [->]  (covn)
      (addm) edge [->]  (covm)
      (addm) edge [->] (b)
      (d)  edge[->] (cfm)
      (b) edge [->] (d)
      (nonn) edge [->]  (cfn)
      (cfn) edge[->] (c)
      (covm) edge [->] (d)
      (covm) edge [->]  (nonsn)
      (nonsn) edge [->]  (nonn);

\draw[cyan,line width=.05cm](addn) edge[->] node[left] {} (e2);

\draw[cyan,line width=.05cm] (covn) edge[<-] node[above right] {} (e2);
\draw[cyan,line width=.05cm] (covn) edge[->] node[midway, above left] {} (b);

\draw[cyan,line width=.05cm] (b)  edge [->] node[ below left] {} (nonm);

\draw[cyan,line width=.05cm,dashed] (nonm) edge[->] (covm);

\node  at (0.25,4) {\cyan{$\lambda^-_2$}};
\node  at (0.75,5.55) {\cyan{$\lambda^-_4$}};
\node  at (0.25,5) {\cyan{$\lambda^-_3$}};
\node  at (1.27,5.2) {\cyan{$\lambda^-_5$}};

\node (aleph1) at (-1,-3) {$\aleph_1$};
\node (addn-f) at (0,-3){$\lambda_1^\bfrak$};
\node (covn-f) at (0,0){$\lambda_3^\bfrak$};
\node (e2-f) at (0.7,-1.5){$\lambda_2^\bfrak$};
\node (nonn-f) at (3.8,-3) {$\lambda_2^\dfrak$} ;
\node (nonsn-f) at (3.2,-3) {$\lambda_3^\dfrak$} ;
\node (cfn-f) at (3.8,0) {$\lambda_1^\dfrak$} ;
\node (addm-f) at (1.4,-3) {$\bullet$} ;
\node (covm-f) at (2.6,-3) {$\lambda_5^\dfrak$} ;
\node (nonm-f) at (1.4,0) {$\lambda_5^\bfrak$} ;
\node (cfm-f) at (2.6,0) {$\bullet$} ;
\node (b-f) at (1.4,-1.5) {$\lambda_4^\bfrak$};
\node (d-f) at (2.6,-1.5) {$\lambda_4^\dfrak$};
\node (c-f) at (5,0) {$\lambda$};

\draw[cyan,line width=.05cm,dashed] (c-f) edge[->] node[above] {$\lambda^-_1$} (addn);

\draw[gray]
   (covn-f) edge [->] (nonm-f)
   (nonm-f)edge [->] (cfm-f)
   (cfm-f)edge [->] (cfn-f)
   (addn-f) edge [->]  (addm-f)
   (addn-f) edge [->]  (covn-f)
   (addm-f) edge [->]  (covm-f)
   (covm-f) edge [->]  (nonsn-f)
   
   (addm-f) edge [->] (b-f)
   (d-f)  edge[->] (cfm-f)
   (b-f) edge [->] (d-f);

\draw[,cyan,line width=.05cm,dashed] 
              (covn-f) edge [->] (b-f)
              (nonm-f) edge [->] (covm-f)
              (d-f) edge [->] (nonsn-f);

\draw[cyan,line width=.05cm](aleph1) edge[->] (addn-f)
(e2-f) edge [->] (covn-f)
(addn-f) edge[->] (e2-f)
(b-f)  edge [->] (nonm-f)
(covm-f) edge [->] (d-f)
(nonsn-f) edge [->]  (nonn-f)
(nonn-f) edge [->]  (cfn-f)
(cfn-f) edge[->] (c-f);

}
\end{tikzpicture}
\caption{Strategy to force Cicho\'n's maximum: we construct a ccc poset~$\Por$ forcing the constellation at the top, and find a $\sigma$-closed model $N$ such that $\Por\cap N$ forces the constellation at the bottom.}\label{fig:cichoncoll}
\end{figure}
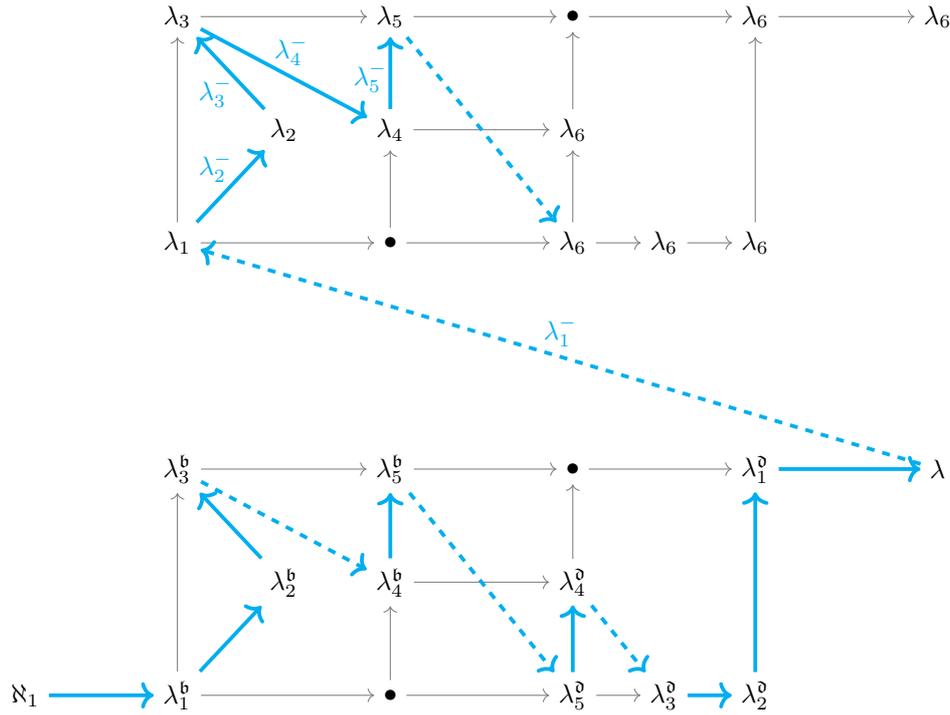

\begin{enumerate}[label = \rm (\alph*)]
\item\label{maxa}
$\cfrak = \lambda_6 = |\Por|$;
\item\label{maxc}
$\Rbf_i \eqT \Cv_{[\lambda_6]^{<\lambda_i}}$ for $1\leq i\leq 4$;
\item\label{maxe}
$\lambda_5\leqT\Cv_\Mwf$, $\lambda_6\leqT \Cv_\Mwf$, and $\Ed\leqT \lambda_6\times\lambda_5$.
\end{enumerate}
For each $1\leq i\leq 4$ let $S_i:=[\lambda_6]^{<\lambda_i}\cap V$, which is a directed preorder. By $\thgch$ above $\lambda^-_1$, we obtain that $S_i\eqT \Cv_{[\lambda_6]^{<\lambda_i}}\cap V$ and that any ccc poset forces $S_i\eqT \Cv_{[\lambda_6]^{<\lambda_i}} \eqT [\lambda_6]^{<\lambda_i}$. Therefore, $\Cv_{[\lambda_6]^{<\lambda_i}}$ can be replaced by $S_i$ in~\ref{maxc}.

On the other hand, under the same dynamics as in~\cite[Sec.~5]{CM22}, we can construct a $\sigma$-closed $N\preceq H_\chi$ of size $\lambda$ such that, for any $1\leq i\leq 4$,
\begin{enumerate}[label = \rm (\alph*$'$)]
\item $S_i\cap N \leqT \prod_{j=i}^5 \lambda^\dfrak_j\times \lambda^\bfrak_j$,
\item $\lambda^\bfrak_i\leqT S_i\cap N$ and $\lambda^\dfrak_i\leqT S_i\cap N$,
\item $\lambda_5\cap N\eqT \lambda^\bfrak_5$ and $\lambda_6\cap N \eqT \lambda^\dfrak_5$.
\end{enumerate}
Then, by~\ref{maxa}--\ref{maxe} and \autoref{KcapN}, $\Por\cap N$ is a required.
\end{varproof}

Just as the proof of prior proof, we can prove~~\autoref{Thm:a1}-\ref{Thm:a2}. 

\section{Open questions}\label{sec:s5}

Regarding on $\lambda$-$\rho$-Knaster property, we ask:

\begin{question}
Does $\lambda$-$\rho$-Knaster property keep the additivity of the null ideal in generic extension? 
\end{question}

In view of~\autoref{e0}, the following is of interest as well:

\begin{question}
Are each one of the following statements consistent with\/ $\thzfc$?
\begin{enumerate}[label=\rm(\arabic*)]
\item $\add(\Nwf)<\efrak^\const<\bfrak$.
\item $\add(\Nwf)<\efrak^\const<\efrak_2^\const<\bfrak$.
\item $\efrak^\const<\efrak_2^\const<\add(\Nwf)<\bfrak$.
\item $\efrak^\const<\add(\Nwf)<\efrak_2^\const<\bfrak$.
\end{enumerate}
\end{question}

One of our original intentions for introducing the property $\Pr^2_{\bar n^*}(\lambda)$ was to solve~\autoref{a0} below, but it did not succeed because random forcing (\autoref{exmp:Pr^2}) satisfies the property $\Pr^2_{\bar n^*}(\lambda)$, so it is still open:

\begin{question}[\cite{kamoeva,kamopred}]\label{a0}
Are each of the following statements consistent with\/ $\thzfc$?
\begin{enumerate}[label=\rm(\arabic*)]
\item\label{a0:a}
$\efrak_2^\const>\cov(\Nwf)$.
\item\label{a0:b}
$\vfa_2^\const<\non(\Nwf)$.
\end{enumerate}
\end{question}

As we mentioned in~\autoref{s0}, it is still open:

\begin{question}[Kada, \cite{Kaunp}]
In\/ $\thzfc$: does\/ $\efrak^\const_2\le\dfrak$ hold?
\end{question}

\subsection*{Acknowledgments}

The authors express their gratitude to Diego A. Mej\'ia for the
productive discussion of the proof of~\autoref{b36}.


{\small
\bibliography{bibli}
\bibliographystyle{alpha}
}


\end{document}